\documentclass[11pt]{article}
\usepackage{times}
\usepackage{amsmath,amsfonts,amstext,latexsym,amssymb,amsbsy,amsopn,amsthm,eucal,slashed,amstext,enumerate,accents}
\usepackage{mathrsfs}

\usepackage{cases}
\usepackage{color}
\usepackage{txfonts}
\usepackage{dsfont}
\usepackage{graphicx}   % for figures
\usepackage{comment}
\usepackage{scalerel,stackengine}
\stackMath
\newcommand\reallywidehat[1]{%
	\savestack{\tmpbox}{\stretchto{%
			\scaleto{%
				\scalerel*[\widthof{\ensuremath{#1}}]{\kern-.6pt\bigwedge\kern-.6pt}%
				{\rule[-\textheight/2]{1ex}{\textheight}}%WIDTH-LIMITED BIG WEDGE
			}{\textheight}% 
		}{0.5ex}}%
	\stackon[1pt]{#1}{\tmpbox}%
}
\allowdisplaybreaks[4]
\usepackage[colorlinks,
linkcolor=blue,      
anchorcolor=blue,  
citecolor=blue,        
]{hyperref}

\usepackage{enumitem,booktabs}

\numberwithin{equation}{section}
\newcommand{\dd}{\mathrm d}
\DeclareMathOperator{\Ric}{Ric}
\DeclareMathOperator{\Div}{div}
\DeclareMathOperator{\tr}{tr}

\newtheorem{theorem}{Theorem}[section]
\newtheorem{lemma}[theorem]{Lemma}
\newtheorem{proposition}[theorem]{Proposition}
\newtheorem{corollary}[theorem]{Corollary}
\theoremstyle{definition}

\theoremstyle{remark}
\newtheorem{remark}[theorem]{Remark}

\begin{document}
	
	\title{\bf Gradient estimates and monotonicity formulae
		for the $p$-filtration equation on Riemannian manifolds}
		
	\author{Yu-Zhao Wang\thanks{The research of Y.-Z. Wang
			was supported by the Fundamental Research Program
			of Shanxi Province (Grant No.~202303021211001).}
		\quad Jian-Hua Hao}

	\date{}
	\maketitle
	
	\begin{abstract}
We derive Li--Yau type estimates on closed Riemannian manifolds
and local Hamilton-type gradient estimates on complete manifolds
for the $p$-filtration equation \(u_t=\Delta_pF(u)\).
Under nonnegative Ricci curvature and a structural condition on $F$, we establish a $\mathcal W$-entropy formula, entropy-power concavity and local entropy rigidity.
We also prove a parabolic frequency monotonicity formula
without curvature assumptions.
The results include homogeneous doubly nonlinear diffusion
and extend to nonhomogeneous power sums and the arctangent law.
\end{abstract}
	
	\vspace{2mm}
	
	\vspace{2mm}

	\section{Introduction and main results}
\label{ms:intro:section}

Gradient estimates and entropy formulas play important roles
in geometric analysis.
In 1986, Li and Yau \cite{LY} proved that every positive solution of the heat equation
\(
\partial_tu=\Delta u
\)
on a complete Riemannian manifold with nonnegative Ricci
curvature satisfies
\begin{equation}\label{ms:intro:LY-heat}
\frac{|\nabla u|^2}{u^2}-\frac{u_t}{u}
\le\frac{n}{2t}.
\end{equation}
The Li--Yau type gradient estimates have been studied for many
other parabolic equations, including sharp Li--Yau type
estimates for several nonlinear diffusion equations;
see, for example, \cite{Hamilton,LNVV,KN,WC}.

In 1994, Yau \cite{Yau94} considered the filtration equation
$\partial_tu=\Delta F(u)$ on a closed Riemannian manifold.
For $F\in C^2((0,\infty))$ with $F'>0$, define
$G$ by $G'(s)=F'(s)/s$.
For a constant $\alpha\ne0$ and a function
$c\in C^1((0,T))$, he obtained the estimate
\begin{equation}\label{ms:intro:Yau}
|\nabla G(u)|^2-\alpha\,\partial_tG(u)\le c(t),
\end{equation}
provided that the Ricci curvature is nonnegative and a quadratic condition
involving $F$, $\alpha$, and $c$ is satisfied.
The filtration equation is also discussed in \cite{VazquezPME}, where
it is called the generalized porous-medium equation.
In 2008, Ma--Zhao--Song \cite{Ma} extended this study
to equations with a reaction term,
$\partial_tu=\Delta F(u)+H(u)$. They obtained differential Harnack estimates under
suitable assumptions on $F$ and $H$, and also established
Hamilton-type gradient estimates for this equation.
More recently, Hao--Wang \cite{HaoWang} used Nash--Moser
iteration to derive local differential Harnack estimates
for filtration equations on complete Riemannian manifolds
with Ricci curvature bounded below.

For $F(s)=s^m$,
Lu--Ni--V\'azquez--Villani \cite{LNVV} obtained
Aronson--B\'enilan and Li--Yau type estimates in 2009.
In particular, on a closed manifold with $\Ric\ge0$, every positive smooth solution of $u_t=\Delta(u^m)$ satisfies
\begin{equation}\label{ms:intro:LNVV-AB}
\Delta v\ge-\frac{n}{[n(m-1)+2]t},
\end{equation}
where
\(
v=\frac{m}{m-1}u^{m-1},
m>\max\left\{0,1-\frac2n\right\},
 m\ne1.
\)

The preceding results concern the case $p=2$ of the
$p$-Laplacian
\(
\Delta_pf=\Div\bigl(|\nabla f|^{p-2}\nabla f\bigr).
\)
In 2009, Kotschwar--Ni \cite{KN} studied the
$p$-Laplacian type parabolic equations
\(
\partial_tu=\Delta_pu
\)
and
\(
\partial_t(u^{p-1})=(p-1)^{p-1}\Delta_pu.
\)
Their results include global Li--Yau type estimates
under nonnegative Ricci curvature.
For the second equation, they also derived local estimates
under a lower bound on the sectional curvature and
additional assumptions, as well as an entropy formula.

These results naturally lead to the study of the
$p$-Laplacian combined with the power function $F(s)=s^m$.
In 2014, Wang--Chen \cite{WC} studied the doubly nonlinear
diffusion equation $\partial_tu=\Delta_p(u^m)$. On a closed manifold with $\Ric\ge0$, for $p>1$ and
$b=m-1/(p-1)>0$, they proved
\begin{equation}\label{DNDE}
\frac{|\nabla v|^p}v-\frac{v_t}v
\le
\frac{nb}
{\bigl[(p-1)nb+p\bigr]t},
\end{equation}
where
\(
v=\frac{m}{b}
u^{b}.
\)
They also obtained an entropy monotonicity formula
for this equation; see \eqref{ms:intro:WC-formula} below.

In 1993, Hamilton \cite{Hamilton} established a matrix
Harnack inequality on closed manifolds with nonnegative
sectional curvature and parallel Ricci curvature.
Taking its trace recovers the Li--Yau estimate.
In the same work, he also proved that every positive solution
of the heat equation on a closed manifold with $\Ric\ge0$
satisfying $u\le A$ obeys
\begin{equation}\label{ms:intro:Hamilton-heat}
t\frac{|\nabla u|^2}{u^2}\le\log\frac Au.
\end{equation}
In 2006, Souplet--Zhang \cite{SZ} obtained a local
Hamilton-type estimate for the heat equation on complete
manifolds with Ricci curvature bounded below.
In 2012, Xu \cite{Xu} considered Hamilton-type estimates
for the filtration equation $u_t=\Delta F(u)$.
Under suitable structural conditions on $F$,
he proved
\begin{equation}\label{ms:intro:Xu}
\sup_{Q_{R/2,T/2}}
\frac{|\nabla G(u)|}{\alpha-G(u)}
\le
C\left(\frac1R+\frac1{\sqrt T}+\sqrt\kappa\right),
\end{equation}
where
\(
Q_{R,T}=B(x_0,R)\times(t_0-T,t_0],
\Ric\ge-\kappa g,
G'(s)=F'(s)/s.
\)
The constant $C$ depends on the dimension and the bounds
in the structural assumptions.
Hamilton-type estimates for other nonlinear parabolic
equations have also been studied; see
\cite{Ma,Wang2026}.

It is natural to consider the $p$-Laplacian combined
with a general increasing function $F$ on a Riemannian
manifold, leading to the $p$-filtration equation
\begin{equation}\label{ms:intro:PF}
	\partial_tu=\Delta_pF(u),\qquad p>1.
\end{equation}
Although this class of equations has a broad range
of applications, its geometric estimates and monotonicity
formulas on Riemannian manifolds remain less studied
for general $F$ and $p\ne2$.
V\'azquez's monographs \cite{VazquezSmoothing,VazquezPME}
discuss models in this class arising in non-Newtonian
fluid mechanics, flows through porous media,
and image processing.

The equation can be written in the mass-conservation form
\begin{equation}\label{ms:intro:mass-conservation}
	\partial_tu+\Div(u\mathbf V)=0,
\end{equation}
where
\(
\mathbf V
=-u^{-1}|\nabla F(u)|^{p-2}\nabla F(u).
\)
Define the pressure by
\begin{equation}\label{ms:intro:pressure}
	v(u):=c_0+\int_{u_*}^{u}
	\frac{F'(s)}{s^{1/(p-1)}}\,\dd s,
\end{equation}
where $u_*$ is a fixed positive reference value and
$c_0\in\mathbb R$.
Then
\begin{equation}\label{ms:intro:pressure-gradient}
	\nabla v=u^{-1/(p-1)}\nabla F(u),
	\qquad
	\nabla F(u)=u^{1/(p-1)}\nabla v.
\end{equation}
Thus $\mathbf V=-|\nabla v|^{p-2}\nabla v$,
and \eqref{ms:intro:mass-conservation} becomes
\begin{equation}\label{ms:intro:continuity}
	\partial_tu
	=\Div\bigl(u|\nabla v|^{p-2}\nabla v\bigr).
\end{equation}
The first part of the paper concerns gradient estimates,
beginning with the following Li--Yau type estimate.

\begin{theorem}[Li--Yau type estimates]
	\label{ms:thm:LY}
	Let $(M^n,g)$ be closed with $\Ric\ge0$, and let $p\ge2$.
	Assume that $F\in C^3$ and $F'>0$ on a neighborhood
	of $\mathcal U=[u_-,u_+]\Subset(0,\infty)$.
	Let $u>0$ be a smooth solution of \eqref{ms:intro:PF} on $M\times(0,T)$
	with values in $\mathcal U$, and set
	\(
	\varphi=\log u,
	\psi(\varphi)=v(e^\varphi),
	\)
	where $v$ is the pressure function defined
	by \eqref{ms:intro:pressure}.
	Fix $\alpha\ge1$ and a nonnegative function
	$c\in C^1((0,T))$.
	Suppose that
	\begin{equation}\label{ms:eq:LY-polynomial}
		\begin{aligned}
			&\frac{1-\alpha}{\alpha}s^{\frac{p-2}{p-1}}
			\left[
			(p-1)sF''(s)
			+\left(\frac{p(\alpha-1)}{n\alpha}+p-2\right)F'(s)
			\right]y^2+
			\left[
			\frac{2p(1-\alpha)}{n\alpha^2}
			-\frac{\alpha+p-2}{\alpha}
			\left(\frac{sF''(s)}{F'(s)}+\frac{p-2}{p-1}\right)
			\right]c(t)y\\
			&\quad-
			\frac{c(t)^2}{\alpha s^{\frac{p-2}{p-1}}F'(s)}
			\left[
			\frac{p}{n\alpha}
			+\frac{sF''(s)}{F'(s)}
			+\frac{p-2}{p-1}
			\right]
			-c'(t)\le0
		\end{aligned}
	\end{equation}
	for every $s\in\mathcal U$, $y\ge0$, and $0<t<T$.
	If, at some $\tau\in(0,T)$,
	\[
	|\nabla\psi|^p-\alpha\psi_t\le c(\tau)
	\qquad\text{on }M\times\{\tau\},
	\]
	then
	\begin{equation}\label{ms:eq:LY-barrier}
		|\nabla\psi|^p-\alpha\psi_t\le c(t)
		\qquad\text{on }M\times[\tau,T).
	\end{equation}
\end{theorem}

\begin{remark}\label{ms:rem:LY-models}
	\leavevmode\par
	\begin{itemize}
		\item
		For $p=2$, one has $\psi=G(u)$ up to an additive constant,
		so \eqref{ms:eq:LY-barrier} is Yau's estimate \eqref{ms:intro:Yau},
		and condition \eqref{ms:eq:LY-polynomial}
		coincides with the condition of Ma--Zhao--Song
		\cite[Theorem~3]{Ma} when $H\equiv0$.
		
		\item
		The criterion also applies to positive sums
		$F(s)=\sum_j c_js^{m_j}$ and to $F(s)=\arctan s$; see
		Section~\ref{ms:examples:section}.
	\end{itemize}
\end{remark}

An Aronson--B\'enilan estimate for $\Delta_pv$ follows under an additional monotonicity condition
on the pressure coefficient.

\begin{theorem}[Aronson--B\'enilan estimate]
	\label{ms:thm:AB}
	Let $(M^n,g)$ be closed with $\Ric\ge0$, and let $p\ge2$.
	Assume that $F\in C^3$ and $F'>0$ on a neighborhood of
	$\mathcal U=[u_-,u_+]\Subset(0,\infty)$.
	Let $u>0$ be a smooth solution of
	\eqref{ms:intro:PF} on $M\times(0,T)$
	with values in $\mathcal U$.
	Let $v$ be the pressure defined by
	\eqref{ms:intro:pressure}, and set
	$\varphi=\log u$ and $\psi(\varphi)=v(e^\varphi)$.
	Define
	\begin{equation}\label{ms:eq:b-def}
		b(s)=\frac{p-2}{p-1}+\frac{sF''(s)}{F'(s)},
	\end{equation}
	so that $b(u)=\psi''(\varphi)/\psi'(\varphi)$.
	Suppose that $b$ is nonincreasing on $\mathcal U$
	and that, for some constant $b_0$,
	\begin{equation}\label{ms:eq:AB-condition}
		b(s)\ge b_0\quad(s\in\mathcal U),
		\qquad
		\varrho:=(p-1)b_0+\frac pn>0.
	\end{equation}
	Then, on $M\times(0,T)$,
	\begin{equation}\label{ms:eq:AB-bound}
		\Delta_p\psi\ge-\frac1{\varrho t},
		\qquad\text{equivalently}\qquad
		\frac{|\nabla \psi|^p}{\psi'}-\frac{\psi_t}{\psi'}
		\le\frac{1}{\varrho t}.
	\end{equation}
\end{theorem}

\begin{remark}\label{ms:rem:AB}
	\leavevmode\par
	\begin{itemize}
		\item
		For $F(s)=s^m$, the coefficient
		$b(s)=b=m-1/(p-1)$ is constant, and the theorem applies
		whenever $m>0, p\ge 2$ and $b>-\frac p{p-1}\frac1n$.
		In particular, for $b>0$, choosing
		$v=\frac mb u^b$
		gives the estimate \eqref{DNDE} of Wang--Chen \cite{WC}.
		For $p=2$, the theorem gives the
		Aronson--B\'enilan estimate \eqref{ms:intro:LNVV-AB} of
		Lu--Ni--V\'azquez--Villani \cite{LNVV}
		in the range $m>\max\{0,1-2/n\}$, $m\ne1$.
		
		\item
		The monotonicity condition also holds for $F(s)=\log s$, where
		$b=-1/(p-1)$ and hence $\Delta_pv\ge-n/((p-n)t)$ for $p\ge2$,
		$p>n$, and for $F(s)=\arctan s$; see
		Section~\ref{ms:examples:section}.
	\end{itemize}
\end{remark}

The next result gives the local Hamilton-type estimate.

\begin{theorem}[Local Hamilton-type estimate]
	\label{ms:thm:Hamilton}
	Let $(M^n,g)$ be complete, let $p>1$, and suppose
	$\operatorname{Sec}\ge-\kappa$ on $B(x_0,R)$,
	where $\kappa\ge0$.
	Assume that $F\in C^3$ and $F'>0$ on a neighborhood
	of $\mathcal U=[u_-,u_+]\Subset(0,\infty)$.
	Let $u>0$ be a smooth solution of \eqref{ms:intro:PF} in
	\(
	Q_{R,T}=B(x_0,R)\times(t_0-T,t_0]
	\)
	with values in $\mathcal U$.
	For the pressure $v$ defined by
	\eqref{ms:intro:pressure}, set
	\(
	\varphi=\log u,
	\psi(\varphi)=v(e^\varphi).
	\)
	Suppose that there exist constants
	$\alpha\in\mathbb R$ and $q,\delta,\beta>0$ such that
	\[
	\alpha-\sup_{s\in\mathcal U}v(s)\ge\delta>0,
	\]
	and
	\begin{equation}\label{ms:eq:Hamilton-condition}
		\begin{aligned}
			&q(p-1)
			+q\left(p-1+\frac{(n-1)(p-2)}2\right)
			\left(
			\frac{sF''(s)}{F'(s)}+\frac{p-2}{p-1}
			\right)\\
			&+
			q\left[
			p-1-\frac qp
			\left((p-1)^2+\frac{(n-1)(p-2)^2}{4}\right)
			\right]
			\frac{s^{\frac{p-2}{p-1}}F'(s)}{\alpha-v(s)}-
			\frac{(n-1)p[\alpha-v(s)]}
			{4s^{\frac{p-2}{p-1}}F'(s)}
			\left(
			\frac{sF''(s)}{F'(s)}+\frac{p-2}{p-1}
			\right)^2
			\ge\beta
		\end{aligned}
	\end{equation}
	for every $s\in\mathcal U$.
	Then
	\begin{equation}\label{ms:eq:Hamilton-bound}
		\sup_{Q_{R/2,T/2}}
		\frac{|\nabla\psi|}{(\alpha-\psi)^{q/p}}
		\le C\left(
		\frac1R+\frac1{T^{1/p}}+\sqrt{\kappa}
		\right).
	\end{equation}
	The constant $C$ depends only on $n,p,q,\delta,\beta$ and an upper
	bound $\Lambda$ for $\sup_{\mathcal U}\psi'$, $\sup_{\mathcal U}|b|$
	and $\operatorname{osc}_{\mathcal U}v$; in particular, it is
	independent of $R$, $T$, $\kappa$ and $x_0$.
\end{theorem}

\begin{remark}\label{ms:rem:Hamilton-models}
	\leavevmode\par
	\begin{itemize}
		\item
		For $p=2$ and $F(s)=s$, one has $\psi=\log u$, and the choices $q=1$ and
		$q=2$ give local versions of Hamilton's estimate
		\eqref{ms:intro:Hamilton-heat} and of the Souplet--Zhang estimate
		\cite{SZ}, respectively.
		For $F(s)=s^{1/(p-1)}$, where $b\equiv0$, condition
		\eqref{ms:eq:Hamilton-condition} holds with $\beta=q(p-1)$
		for every $q$ satisfying \eqref{ms:examples:Hamilton-q}.
		
		\item
		For $p=q=2$, the quotient in \eqref{ms:eq:Hamilton-bound}
		is the one in Xu's estimate \eqref{ms:intro:Xu}, and condition
		\eqref{ms:eq:Hamilton-condition} reduces to
		$2+2b-(n-1)b^2/(2d)\ge\beta>0$, where $b=sF''/F'$ and
		$d=F'/(\alpha-v)$, whereas condition~{\rm(C)} in
		\cite[Theorem~1.1]{Xu} requires $2+2b-(n-1)b^2/d\ge\gamma>0$.
		For a fixed pressure gap, our condition is therefore weaker than
		Xu's stated condition, while retaining the coefficient already
		obtained in \cite[equation~(15)]{Xu}.
		
	\end{itemize}
\end{remark}

Entropy formulas provide integral estimates that complement
the preceding gradient estimates. In 2002, in his work on the Ricci
flow, Perelman \cite{Perelman} introduced the $\mathcal W$-entropy and
proved its monotonicity. Following this idea, in 2009,
Lu--Ni--V\'azquez--Villani \cite{LNVV} established
$\mathcal W$-entropy formulas for the porous-medium and
fast-diffusion equations, and Kotschwar--Ni \cite{KN} obtained an
entropy formula for the $p$-heat equation. In 2012, X.-D. Li  \cite{LiWitten} proved the $W$-entropy formula for the Witten Laplacian heat equation via Bakry--\'Emery Ricci curvature. In 2014, Wang--Chen
\cite{WC} combined these two directions for the doubly nonlinear
equation $\partial_tu=\Delta_p(u^m)$, and their formula is the
starting point of the present work. Let $(M^n,g)$ be closed with $\Ric\ge0$,
let $p>1$ and $b=m-1/(p-1)>0$, and let
\[
v=\frac mbu^b,\qquad a=\frac{nb}{n(p-1)b+p}
\]
be the pressure and the constant in \eqref{DNDE}. Wang--Chen defined
the Nash-type entropy and the Perelman-type entropy by
\begin{equation}\label{ms:intro:WC-entropy}
	\mathcal N_u(t)=-t^a\int_Mvu\,\dd V,
	\qquad
	\mathcal W_u(t)=t\frac{\dd}{\dd t}\mathcal N_u(t)+\mathcal N_u(t)
	=t^{a+1}\int_M\left((b+1)\frac{|\nabla v|^p}{v}-\frac{a+1}{t}\right)vu\,\dd V,
\end{equation}
and proved that
\begin{equation}\label{ms:intro:WC-formula}
	\begin{aligned}
		\frac{\dd}{\dd t}\mathcal W_u
		={}&-pbt^{a+1}\int_M\left(
		\left||\nabla v|^{p-2}\nabla_i\nabla_jv+\frac{a}{nbt}a_{ij}\right|_A^2
		+|\nabla v|^{2p-4}\Ric(\nabla v,\nabla v)\right)vu\,\dd V\\
		&-pt^{a+1}\int_M\left(b\Delta_pv+\frac at\right)^2vu\,\dd V\le0,
	\end{aligned}
\end{equation}
where $A^{ij}=g^{ij}+(p-2)v^iv^j/|\nabla v|^2$ is the linearization
matrix of $\Delta_p$, $(a_{ij})$ is its inverse, and $|\cdot|_A$ is the
induced norm, as in \cite{WC}. For $p=2$,
\eqref{ms:intro:WC-formula} reduces to the porous-medium formula of
\cite{LNVV}, and as $b\to0$ it becomes the formula of \cite{KN}.

The second derivative of an entropy also leads to
entropy-power concavity.
In 1985, Costa \cite{Costa} proved that the Shannon entropy power is
concave along the Euclidean heat flow, and in 2014,
Savar\'e--Toscani \cite{ST} established the corresponding result for the R\'enyi entropy power:
probability solutions of $w_t=\Delta(w^m)$
with $m>1$ and finite R\'enyi entropy satisfy
\begin{equation}\label{ms:intro:Renyi}
	\mathcal N_m(t)
	=\left(\int_{\mathbb R^n}w^m\,\dd x\right)^{
		-\frac{m-1+2/n}{m-1}},
	\qquad
	\mathcal N_m''(t)\le0.
\end{equation}
S. Li and X.-D. Li \cite{LL,LLRenyi} later studied Shannon and R\'enyi
entropy powers on Riemannian manifolds, and since 2019 the first
author and his coauthors \cite{WZ,WW,Wang2026} have obtained
$p$-Shannon entropy-power concavity for the $p$-heat equation and
$p$-R\'enyi entropy-power concavity for the doubly
nonlinear diffusion equation on Riemannian manifolds.

For a general increasing function $F$, the scaling structure behind
\eqref{ms:intro:WC-formula} is lost: the equation has no
Barenblatt-type self-similar solutions, and the coefficients in
\eqref{ms:intro:WC-formula} are no longer tied to $F$. We show that a
single scalar quantity $D_a$, defined in \eqref{ms:intro:scalar-defect}
below, the condition \(D_a\ge0\) is sufficient for both
$\mathcal W$-entropy monotonicity and the entropy-power concavity
for \eqref{ms:intro:PF}.

For the entropy results, assume that
$F\in C^3((0,\infty))$, $F'>0$, and
$\int_0^s r^{-1/(p-1)}F'(r)\,\dd r<\infty$ for $s>0$.
Normalize the pressure in \eqref{ms:intro:pressure}
by $v(0)=0$, and set
\begin{equation}\label{ms:intro:energy-notation}
	\begin{aligned}
		U(s)&=\int_0^s v(r)\,\dd r,
		\qquad
		P(s)=sv(s)-U(s)
		=\int_0^s r^{\frac{p-2}{p-1}}F'(r)\,\dd r,\\
		\mathcal E(t)&=\int_M U(u)\,\dd V,
		\qquad
		\mathcal I(t)=\int_M u|\nabla v|^p\,\dd V.
	\end{aligned}
\end{equation}
Then $U(s),P(s)>0$ for $s>0$, and
\eqref{ms:intro:continuity} gives
$\mathcal E'=-\mathcal I$ on a closed manifold.
Following \eqref{ms:intro:WC-entropy}, for $a>0$ we define the
$\mathcal W$-entropy
\begin{align}
	\mathcal W_{F,p,a}(t)
	:=t\frac{\dd}{\dd t}\bigl(-t^a\mathcal E(t)\bigr)-t^a\mathcal E(t)
	=t^{a+1}\mathcal I(t)
	-(a+1)t^a\mathcal E(t),
	\label{ms:intro:W-definition}
\end{align}	
and the entropy power	
\begin{align}
	\mathcal N_{F,p,a}(t)
	:=\mathcal E(t)^{-1/a}.
	\label{ms:intro:N-definition}
\end{align}
For $F(s)=s^m$ with $m>1/(p-1)$,
the corresponding choice of $a$ gives the
R\'enyi entropy power up to a positive constant.

The structural condition below involves the quantity
$sP'(s)-(1-\frac1n)P(s)$, which also appears in McCann's criterion for
displacement convexity of internal energies
\cite{McCann}.
For a general $F$, an additional scalar term
enters the entropy dissipation.
Together, these terms give a condition for both
$\mathcal W$-entropy monotonicity and
entropy-power concavity.
This condition is verified for $p=2$, for positive sums of powers
and for the arctangent law in
Section~\ref{ms:examples:section}.
The first result is the following
$\mathcal W$-entropy formula.

\begin{theorem}[$\mathcal W$-entropy formula]
	\label{ms:thm:entropy}
		Let $(M^n,g)$ be closed with $\Ric\ge0$, and let $p>1$.
	Let $u>0$ be a smooth solution of
	\eqref{ms:intro:PF} on $M\times(0,T)$
	with values in $\mathcal U$, whose pressure flux
	$|\nabla v|^{p-2}\nabla v$ is $C^1$ in space. Use the normalization
	and notation in \eqref{ms:intro:energy-notation} and
	\eqref{ms:intro:WC-formula}.
	For fixed $a>0$, set
	\begin{equation}\label{ms:intro:scalar-defect}
		D_a(s)
		:=sP'(s)-\left(1-\frac1n\right)P(s)
		-\frac{a+1}{pa}\frac{P(s)^2}{U(s)}.
	\end{equation}
	Then the functional $\mathcal W_{F,p,a}$ defined
	in \eqref{ms:intro:W-definition} satisfies
	\begin{equation}\label{ms:ent:Wsquare}
		\begin{aligned}
			\frac{\dd}{\dd t}\mathcal W_{F,p,a}
			={}&-pt^{a+1}\int_M
			P(u)\left(
			\left||\nabla v|^{p-2}\nabla_i\nabla_jv+\frac{aU(u)}{ntP(u)}a_{ij}\right|_A^2
			+|\nabla v|^{2p-4}\Ric(\nabla v,\nabla v)
			\right)\dd V\\
			&-pt^{a+1}\int_M
			\bigl[uP'(u)-P(u)-D_a(u)\bigr]
			\left(\Delta_pv+\frac{aU(u)}{tP(u)}\right)^2\dd V
			-pt^{a+1}\int_MD_a(u)(\Delta_pv)^2\,\dd V.
		\end{aligned}
	\end{equation}
	If, in addition,
	\begin{equation}\label{ms:ent:scalar-condition}
		D_a(s)\ge0,\qquad s\in\mathcal U,
	\end{equation}
	then $\mathcal W_{F,p,a}'(t)<0$ for every $t\in(0,T)$.
\end{theorem}

Identity \eqref{ms:ent:Wsquare} has the same form as
\eqref{ms:intro:WC-formula}. Since
$uP'-P-D_a=\frac{a+1}{pa}\frac{P^2}{U}-\frac Pn$, this coefficient can be
negative even when $D_a>0$; for $F(s)=s^m$ this happens when $nb<p$ and
$a>nb/(p-nb)$. Thus the sign of $\mathcal W_{F,p,a}'$ is not read off
term by term from \eqref{ms:ent:Wsquare}. Instead, by
\eqref{ms:ent:W-dissipation},
\[
\begin{aligned}
\frac{\dd}{\dd t}\mathcal W_{F,p,a}
={}&-pt^{a+1}\int_M\Big\{P(u)\Big(\Big||\nabla v|^{p-2}\nabla_i\nabla_jv
-\frac{\Delta_pv}{n}a_{ij}\Big|_A^2+|\nabla v|^{2p-4}\Ric(\nabla v,\nabla v)\Big)\\
&\qquad+D_a(u)(\Delta_pv)^2+\frac{a+1}{pa}\frac{P(u)^2}{U(u)}
\Big(\Delta_pv+\frac{aU(u)}{tP(u)}\Big)^2\Big\}\,\dd V,
\end{aligned}
\]
in which every term of the integrand is nonnegative under
\eqref{ms:ent:scalar-condition}.

The same computation also yields
entropy-power concavity under the same structural condition.

\begin{theorem}[Concavity of entropy power]
	\label{ms:thm:entropy-power}
	Under the hypotheses of Theorem~\ref{ms:thm:entropy}
	and \eqref{ms:ent:scalar-condition},
	the entropy power $\mathcal N_{F,p,a}$ defined in
	\eqref{ms:intro:N-definition} satisfies
	\[
	\mathcal N_{F,p,a}'
	=\frac1a\mathcal E^{-1/a-1}\mathcal I\ge0
	\]
	and
	\begin{equation}\label{ms:ent:power-second}
		\begin{aligned}
			\mathcal N_{F,p,a}''
			={}&-\frac pa\mathcal E^{-1/a-1}
			\int_M
			\Bigg\{
			P(u)\left[
			\left||\nabla v|^{p-2}\nabla_i\nabla_jv-\frac{\Delta_pv}{n}a_{ij}\right|_A^2
			+|\nabla v|^{2p-4}\Ric(\nabla v,\nabla v)
			\right]\\
			&\qquad+D_a(u)(\Delta_pv)^2
			+\frac{a+1}{pa}U(u)
			\left(
			\frac{P(u)\Delta_pv}{U(u)}
			+\frac{\mathcal I}{\mathcal E}
			\right)^2
			\Bigg\}\,\dd V
			\le0.
		\end{aligned}
	\end{equation}
	At any $t_0\in(0,T)$, equality
	$\mathcal N_{F,p,a}''(t_0)=0$ holds if and only if
	$u(\cdot,t_0)$ is constant on each connected component
	of $M$. In this case, the solution remains stationary
	for $t\ge t_0$.
\end{theorem}

\begin{remark}\label{ms:rem:entropy-models}
	Let $F(s)=s^m$ with $b=m-1/(p-1)>0$, and let
	$a=nb/(n(p-1)b+p)$. Then $v=\frac mbs^b$, $U=sv/(b+1)$, $P=bU$,
	$sP'-P=bP$ and $D_a\equiv0$; hence $aU/(nP)=a/(nb)$ and
	$[uP'-P](\Delta_pv+\frac{aU}{tP})^2=U(b\Delta_pv+\frac at)^2$. Since
	the energy density in \eqref{ms:intro:WC-entropy} is $uv=(b+1)U(u)$,
	one has $\mathcal N_u=-(b+1)t^a\mathcal E$ and
	$\mathcal W_u=(b+1)\mathcal W_{F,p,a}$, and \eqref{ms:ent:Wsquare}
	becomes term by term the formula \eqref{ms:intro:WC-formula} of
	Wang--Chen~\cite{WC}. Moreover,
	$\mathcal N_{F,p,a}=(b(b+1)/m)^{1/a}
	(\int_Mu^{b+1}\,\dd V)^{-(p-1)-p/(nb)}$, so
	Theorem~\ref{ms:thm:entropy-power} recovers the $p$-R\'enyi
	entropy-power concavity of Wang--Wang~\cite{WW}. For $p=2$ and $m>1$,
	these specialize to the porous-medium formula of~\cite{LNVV} and, up
	to a positive constant, to the R\'enyi entropy power
	\eqref{ms:intro:Renyi} of Savar\'e--Toscani~\cite{ST} on $M$, whose
	concavity on closed manifolds with $\Ric\ge0$ is also established
	in~\cite{LLRenyi}.
\end{remark}

The preceding theorems characterize the equality cases
on closed manifolds. We next examine the simultaneous
vanishing of the entropy-production terms on an open set.
The resulting rigidity determines both the diffusion law
on the attained density interval and the local geometry.
The latter is related to the classical structure of
homothetic gradient fields \cite{Tashiro}.

\begin{theorem}[Local entropy rigidity]
	\label{ms:thm:rigidity}
	Let $n\ge2$, $p>1$, and $a>0$.
	Let $\Omega$ be a connected open subset of a
	Riemannian manifold with $\Ric\ge0$.
	Assume that $F$, $U$, and $P$ satisfy the constitutive
	assumptions of Theorem~\ref{ms:thm:entropy},
	and let $u>0$ be smooth on $\Omega$, with
	\(
	|\nabla v|>0.
	\)
	Suppose that, for a constant $\lambda>0$,
	\begin{equation}\label{ms:intro:local-equality-system}
		\begin{gathered}
			|\nabla v|^{2p-4}|\nabla\nabla v|_A^2
			-\frac{(\Delta_pv)^2}{n}=0,
			\qquad |\nabla v|^{2p-4}\Ric(\nabla v,\nabla v)=0,\\
			D_a(u)=0,
			\qquad
			\frac{P(u)\Delta_pv}{U(u)}+\lambda=0
		\end{gathered}
	\end{equation}
	throughout $\Omega$.
	Then
	\[
	0<a<\frac1{p-1},
	\qquad
	\delta=\frac{pa}{n[1-(p-1)a]}>0.
	\]
	There exist constants $C_U,C_F>0$ and $C_0\in\mathbb R$
	such that, on the attained density interval $u(\Omega)$,
	\begin{equation}\label{ms:intro:rigidity-power-laws}
		U(s)=C_U s^{1+\delta},
		\qquad
		F(s)=C_Fs^{\delta+1/(p-1)}+C_0.
	\end{equation}
	Moreover,
	\begin{equation}\label{ms:intro:local-cone}
		\begin{gathered}
			\nabla\bigl(|\nabla v|^{p-2}\nabla v\bigr)
			=-\frac{\lambda}{n\delta}\operatorname{id},\qquad
			g=\dd r^2+r^2g_\Sigma
			\quad\text{locally on }\Omega,
			\qquad
			r=\frac{n\delta}{\lambda}|\nabla v|^{p-1},
		\end{gathered}
	\end{equation}
	where $g_\Sigma$ is a metric on a local level
	hypersurface of $r$.
\end{theorem}

The final part of the paper concerns the time dependence of
 a solution relative to a constant state. Parabolic frequency relates
the rate of energy decay to the energy itself, leading to temporal
interpolation and nonvanishing results.

In 1967, Agmon--Nirenberg \cite{AN} developed logarithmic-convexity
methods to obtain lower bounds and prove uniqueness for solutions
of evolution equations in Hilbert spaces. For the heat equation $w_t=\Delta w$ on closed
manifold, this principle has a simple energy interpretation. Set
\[
H(t)=\int_M w^2\,\dd V,\qquad
E(t)=\int_M|\nabla w|^2\,\dd V,\qquad
N_{\mathrm h}(t)=\frac{E(t)}{H(t)}
\]
on an interval where $H>0$. Integration by parts gives
$H'=-2E$ and $E'=-2\int_M(\Delta w)^2\,\dd V$, hence
\begin{equation}\label{ms:intro:heat-frequency}
N_{\mathrm h}'
=-\frac{2}{H}\int_M(\Delta w+N_{\mathrm h}w)^2\,\dd V\le0,
\qquad
(\log H)''=-2N_{\mathrm h}'\ge0.
\end{equation}
The negative Dirichlet-energy convention used below reverses the
direction of frequency monotonicity.

In 1996, Poon \cite{Poon} developed parabolic frequency methods
for unique continuation using localized space--time quantities.
In 2022, Colding--Minicozzi \cite{CM} established frequency
monotonicity for drift heat equations on manifolds without curvature
assumptions. In 2026, Sun--S\"urig \cite{SS} obtained frequency
formulas for $u_t=\mathcal L_{p,\varphi}(u^q)$ on weighted complete
manifolds, including sign-changing solutions and controlled
perturbations. Also in 2026, Hao--Wang \cite{HaoWang} treated general
filtration laws in $u_t=\Delta F(u)$ through an energy determined
by a primitive of $F(u)-F(\ell)$.
For $u_t=\Delta_pF(u)$, this primitive is combined with the
$p$-Dirichlet dissipation. A scalar condition on the diffusion law
then yields frequency monotonicity and time convexity without a
curvature hypothesis, while an exact square decomposition records
the dissipation terms.

\begin{theorem}[Parabolic frequency monotonicity]
	\label{ms:thm:frequency}
	Let $(M,g)$ be closed, let $p>1$, and let $u>0$
	be a smooth solution of \eqref{ms:intro:PF}
	on $M\times(0,T)$ with values in $\mathcal U$.
	Fix a reference density $\ell$.
	Assume that $F$ is continuous and strictly increasing
	on an interval containing $\ell$ and $\mathcal U$,
	and that $F\in C^2$ with $F'>0$ on a neighborhood
	of $\mathcal U$.
	Set
	\begin{equation}\label{ms:intro:frequency-energy}
		\begin{gathered}
			\Psi_\ell(s)
			=\int_\ell^s[F(r)-F(\ell)]\,\dd r,
			I(t)=2\int_M\Psi_\ell(u)\,\dd V,
			\qquad
			D(t)=-\int_M|\nabla F(u)|^p\,\dd V.
		\end{gathered}
	\end{equation}
	Suppose that, for some $\kappa>0$,
	\begin{equation}\label{ms:freq:structure}
		[F(s)-F(\ell)]^2
		\le\kappa F'(s)\Psi_\ell(s),
		\qquad s\in\mathcal U.
	\end{equation}
	Then, for every $0<\gamma\le p/\kappa$,
	the frequency
	\begin{equation}\label{ms:freq:frequency}
		N_\gamma(t)=\frac{D(t)}{I(t)^\gamma}
	\end{equation}
	is nondecreasing on every time interval on which
	$I>0$. More precisely,
	\begin{equation}\label{ms:freq:main-derivative}
		\begin{aligned}
			N_\gamma'
			={}&pI^{-\gamma}
			\int_MF'(u)
			\left[
			u_t-\frac{2D}{\kappa I}
			\frac{F(u)-F(\ell)}{F'(u)}
			\right]^2\,\dd V\\
			&+\frac{4pD^2}{\kappa^2I^{\gamma+2}}
			\int_M
			\left[
			\kappa\Psi_\ell(u)
			-\frac{[F(u)-F(\ell)]^2}{F'(u)}
			\right]\,\dd V+2\left(\frac p\kappa-\gamma\right)
			\frac{D^2}{I^{\gamma+1}}
			\ge0.
		\end{aligned}
	\end{equation}
\end{theorem}

\begin{remark}\label{ms:rem:frequency-models}
	\leavevmode\par
	\begin{itemize}
		\item
		For $F(s)=s^m$ with $m>0$ and $\ell=0$, one has
		$\Psi_0(s)=s^{m+1}/(m+1)$ and $\kappa=(m+1)/m$, so
		Theorem~\ref{ms:thm:frequency} with $\gamma=pm/(m+1)$ gives the
		monotonicity of
		$-\int_M|\nabla u^m|^p\,\dd V\,\bigl(\int_Mu^{m+1}\,\dd V\bigr)^{-pm/(m+1)}$.
		This is the unweighted closed-manifold case of the generalized
		frequency formula of Sun--S\"urig~\cite[Lemma~2.1]{SS}.

		\item
		For $p=2$ and general $F$, the theorem recovers the
		entropy--dissipation frequency monotonicity of
		Hao--Wang~\cite[Theorem~1.6]{HaoWang}, with the opposite
		sign convention and a constant normalization factor.
		
		\item
		For $p=2$ and $F(s)=s$, one has
		$\Psi_\ell(s)=(s-\ell)^2/2$ and $\kappa=2$.
		With $w=u-\ell$ and $\gamma=1$, one has $I=H$ and
		$N_1=-N_{\mathrm h}$, and \eqref{ms:freq:main-derivative} reduces to
		\eqref{ms:intro:heat-frequency}, the constant-weight case
		of Colding--Minicozzi~\cite{CM}.
	\end{itemize}
\end{remark}

For admissible $\gamma$, this formula gives the convexity of $\log I$
when $\gamma=1$ and of $I^{1-\gamma}/(1-\gamma)$ when $\gamma\ne1$; if
$\gamma\ge1$ is admissible and $I(s)>0$, then $I>0$ on $[s,T)$, so a
solution different from the reference constant cannot reach it in
finite time (Section~\ref{ms:freq:section}).

Section~\ref{ms:key:section} introduces the pressure variables and
the linearized $p$-Bochner identities. Sections~\ref{ms:grad:section},
\ref{ms:ent:section} and \ref{ms:freq:section} prove the gradient
estimates, the entropy results and the frequency monotonicity,
respectively, and Section~\ref{ms:examples:section} treats the case
$p=2$, positive sums of powers and the arctangent law.

\section{Preliminaries}
	\label{ms:key:section}
	
	\subsection{The two pressure variables}

		Equation \eqref{ms:intro:PF} is a doubly nonlinear parabolic
	equation. For the existence of weak solutions, see
	\cite{AltLuckhaus,CailletSantambrogio}, and for local boundedness
	and H\"older regularity of degenerate and singular equations, see
	\cite{DiBenedetto,DGV}.

		We carry out the computations on $\mathcal N=\{|\nabla v|>0\}$.
	The pressure transformation \eqref{ms:intro:pressure} yields
	the following two equivalent forms of the pressure equation.
	
	\begin{lemma}\label{ms:lem:pressure}
		On $\mathcal N$, the pressure $v$ satisfies
		\begin{equation}\label{ms:eq:pressure-pde}
			v_t=\bigl(u^{(p-2)/(p-1)}F'(u)\bigr)\Delta_pv+|\nabla v|^p
			=(uv')\Delta_pv+|\nabla v|^p.
		\end{equation}
	\end{lemma}
	
	\begin{proof}
		Since $\nabla v=u^{-1/(p-1)}\nabla F(u)$, one has
		$|\nabla F(u)|^{p-2}\nabla F(u)=u|\nabla v|^{p-2}\nabla v$.
		Consequently
		\begin{align*}
			v_t
			&=F'(u)u^{-1/(p-1)}\Delta_pF(u)\\
			&=F'(u)u^{-1/(p-1)}
			\Div\bigl(u|\nabla v|^{p-2}\nabla v\bigr)\\
			&=F'(u)u^{(p-2)/(p-1)}\Delta_pv
			+F'(u)u^{-1/(p-1)}|\nabla v|^{p-2}
			\langle\nabla u,\nabla v\rangle\\
			&=(uv')\Delta_pv+|\nabla v|^p.
		\end{align*}
	\end{proof}

	\begin{lemma}\label{ms:lem:log-pressure}
		Set $\varphi=\log u$ and $\psi(\varphi)=v(e^\varphi)=v(u)$.
		Then
		\begin{equation}\label{ms:eq:psi-equation}
			\psi_t=\psi'\Delta_p\psi+|\nabla\psi|^p,
		\end{equation}
		where the primes on $\psi$ denote logarithmic derivatives, and
		\begin{equation}\label{ms:eq:psi-derivatives}
			\begin{aligned}
				\psi'&=u^{(p-2)/(p-1)}F'(u),\\
				\psi''&=u^{(2p-3)/(p-1)}F''(u)
				+\frac{p-2}{p-1}u^{(p-2)/(p-1)}F'(u).
			\end{aligned}
		\end{equation}
	\end{lemma}
	
	\begin{proof}
		Writing $\nabla v=\psi'\nabla\varphi$ and using
		$|\nabla F(u)|^{p-2}\nabla F(u)=u|\nabla v|^{p-2}\nabla v$,
		we obtain
		\begin{align}
			\varphi_t
			&=\frac1u\Div\bigl(u|\nabla v|^{p-2}\nabla v\bigr)=\Delta_pv+|\nabla v|^{p-2}
			\langle\nabla v,\nabla\varphi\rangle.
			\label{ms:eq:varphi-equation}
		\end{align}
		Since $\psi_t=\psi'\varphi_t$ and $\psi(x,t)=v(x,t)$,
		multiplication by $\psi'$ gives
		\begin{align*}
			\psi_t
			&=\psi'\Delta_p\psi
			+\psi'|\nabla\psi|^{p-2}
			\langle\nabla\psi,\nabla\varphi\rangle=\psi'\Delta_p\psi+|\nabla\psi|^p.
		\end{align*}
		Finally, $\dd/\dd\varphi=u\,\dd/\dd u$ gives
		\begin{align*}
			\psi'&=uv'(u)=u^{(p-2)/(p-1)}F'(u),\\
			\psi''&=u\frac{\dd}{\dd u}
			\bigl(u^{(p-2)/(p-1)}F'(u)\bigr)=u^{(2p-3)/(p-1)}F''(u)
			+\frac{p-2}{p-1}u^{(p-2)/(p-1)}F'(u).
		\end{align*}
	\end{proof}

	\subsection{Linearization and the \texorpdfstring{$p$}{p}-Bochner formulas}
	
	Define the linearized operator of the $p$-Laplacian at $v$
	\cite{KN,WC} by
	\begin{equation}\label{ms:eq:Lp}
		\mathcal L_p\xi=\Div\bigl(|\nabla v|^{p-2}A(\nabla\xi)\bigr),
		\qquad
		A=\operatorname{id}+(p-2)
		\frac{\nabla v\otimes\nabla v}{|\nabla v|^2},
	\end{equation}
	and its parabolic operator by
	\begin{equation}\label{ms:eq:square}
		\square=\partial_t-(uv')\mathcal L_p.
	\end{equation}
	The coefficient $uv'$ lies outside the divergence. The eigenvalues
	of $A$ are $p-1,1,\ldots,1$, and $A\nabla v=(p-1)\nabla v$.
	For $p=2$, $A=\operatorname{id}$ and $\mathcal L_p=\Delta$;
	terms with a factor $p-2$ are omitted before evaluation at a
	critical point.
	
	\begin{lemma}\label{ms:lem:Lp-properties}
		For functions $f,g$ with the indicated derivatives,
		\begin{align}
			\mathcal L_p(fg)&=(\mathcal L_pf)g+f\mathcal L_pg
			+2|\nabla v|^{p-2}\langle\nabla f,\nabla g\rangle_A,
			\label{ms:eq:Lp-product}\\
			\int_M(\mathcal L_pf)g\,\dd V
			&=\int_Mf(\mathcal L_pg)\,\dd V
			=-\int_M|\nabla v|^{p-2}
			\langle\nabla f,\nabla g\rangle_A\,\dd V,
			\label{ms:eq:Lp-symmetric}
		\end{align}
		where $\langle X,Y\rangle_A=\langle A(X),Y\rangle$.
		In \eqref{ms:eq:Lp-symmetric}, $f$ and $g$ are smooth
		and compactly supported in
		$\{x\in M:|\nabla v(x,t)|>0\}$ at a fixed time $t$.
	\end{lemma}

	\begin{lemma}\label{ms:lem:time-linearization}
		On $\mathcal N$, we have
		\begin{equation}\label{ms:eq:time-linearization}
			\partial_t(\Delta_pv)=\mathcal L_p(v_t).
		\end{equation}
	\end{lemma}
	
	\begin{proof}
		Differentiating the flux gives
		\begin{align*}
			\partial_t(\Delta_pv)
			&=\Div\left(\partial_t(|\nabla v|^{p-2}\nabla v)\right)\\
			&=\Div\left((p-2)|\nabla v|^{p-4}
			\langle\nabla v,\nabla v_t\rangle\nabla v
			+|\nabla v|^{p-2}\nabla v_t\right)\\
			&=\Div\bigl(|\nabla v|^{p-2}A(\nabla v_t)\bigr)
			=\mathcal L_p(v_t).
		\end{align*}
	\end{proof}
	
	\begin{lemma}[$p$-Bochner formulas]\label{ms:lem:bochner}
		Set $w=|\nabla v|^2$. On $\mathcal N$, we have (see \cite{KN,WC})
		\begin{align}
			\mathcal L_pw
			&=2w^{p/2-1}\bigl(|\nabla\nabla v|^2+\Ric(\nabla v,\nabla v)\bigr)
			+2\langle\nabla v,\nabla\Delta_pv\rangle+\left(\frac p2-1\right)w^{p/2-2}|\nabla w|^2,
			\label{ms:eq:bochner-w}\\
			\mathcal L_p(w^{p/2})
			&=pw^{p-2}\bigl(|\nabla\nabla v|_A^2+\Ric(\nabla v,\nabla v)\bigr)
			+pw^{p/2-1}\langle\nabla v,\nabla\Delta_pv\rangle,
			\label{ms:eq:bochner-power}
		\end{align}
		where
		\begin{equation}\label{ms:eq:Hessian-A}
			\begin{aligned}
				|\nabla\nabla v|_A^2
				&:=A^{ik}A^{j\ell}v_{ij}v_{k\ell}=|\nabla\nabla v|^2+\frac{p-2}{2w}|\nabla w|^2
				+\frac{(p-2)^2}{4w^2}\langle\nabla v,\nabla w\rangle^2.
			\end{aligned}
		\end{equation}
	\end{lemma}

	By $A\nabla v=(p-1)\nabla v$ and \eqref{ms:eq:pressure-pde},
	\begin{equation}\label{ms:eq:square-v}
		\mathcal L_pv=(p-1)\Delta_pv,\qquad
		\square v=|\nabla v|^p-(p-2)(uv')\Delta_pv.
	\end{equation}
	We next apply the $p$-Bochner formula to derive the following
	parabolic identities.
	
	\begin{lemma}\label{ms:lem:evolution}
		On $\mathcal N$, we have
		\begin{align}
			\square v_t
			&=pw^{p/2-1}\langle\nabla v,\nabla v_t\rangle
			+(uv')_t\Delta_pv,
			\label{ms:eq:evol-vt}\\
			\square w
			&=2\langle\nabla v,\nabla v_t\rangle
			-2(uv')w^{p/2-1}\Bigl(|\nabla\nabla v|^2
			+\Ric(\nabla v,\nabla v)
			+w^{1-p/2}\langle\nabla v,\nabla\Delta_pv\rangle\Bigr)\notag\\
			&-\left(\frac p2-1\right)(uv')w^{p/2-2}|\nabla w|^2,
			\label{ms:eq:evol-w}\\
			\square(w^{p/2})
			&=pw^{p/2-1}\langle\nabla v,\nabla v_t\rangle
			-p(uv')w^{p-2}\Bigl(|\nabla\nabla v|_A^2
			+\Ric(\nabla v,\nabla v)
			+w^{1-p/2}\langle\nabla v,\nabla\Delta_pv\rangle\Bigr).
			\label{ms:eq:evol-power}
		\end{align}
		
	\end{lemma}
	
	\begin{proof}
		Using \eqref{ms:eq:pressure-pde} and
		\eqref{ms:eq:time-linearization},
		\begin{align*}
			\square v_t
			&=v_{tt}-(uv')\mathcal L_pv_t\\
			&=v_{tt}-\partial_t\bigl((uv')\Delta_pv\bigr)
			+(uv')_t\Delta_pv\\
			&=pw^{p/2-1}\langle\nabla v,\nabla v_t\rangle
			+(uv')_t\Delta_pv.
		\end{align*}
		For \eqref{ms:eq:evol-w}, subtract $(uv')\mathcal L_pw$ from
		$w_t=2\langle\nabla v,\nabla v_t\rangle$ and use
		\eqref{ms:eq:bochner-w}. Substitution of the pressure equation
		also yields its simplified form:
		\begin{align}
			\square w
			&=2\left\langle\nabla v,
			\nabla\bigl((uv')\Delta_pv+w^{p/2}\bigr)\right\rangle
			-2(uv')w^{p/2-1}
			\bigl(|\nabla\nabla v|^2+\Ric(\nabla v,\nabla v)\bigr)
			\notag\\
			&\quad-(uv')\left[2\langle\nabla v,\nabla\Delta_pv\rangle
			+\left(\frac p2-1\right)w^{p/2-2}|\nabla w|^2\right]
			\notag\\
			&=2\langle\nabla v,\nabla(uv')\rangle\Delta_pv
			+pw^{p/2-1}\langle\nabla v,\nabla w\rangle
			-\left(\frac p2-1\right)(uv')w^{p/2-2}|\nabla w|^2
			\notag\\
			&\quad-2(uv')w^{p/2-1}
			\bigl(|\nabla\nabla v|^2+\Ric(\nabla v,\nabla v)\bigr).
			\label{ms:eq:evol-w-simplified}
		\end{align}
		Finally, the chain rule and this expression give
		\begin{align*}
			\square(w^{p/2})
			&=\frac p2w^{p/2-1}\square w
			-\frac p2\left(\frac p2-1\right)(uv')w^{p-3}
			\langle\nabla w,A\nabla w\rangle\\
			&=pw^{p/2-1}\langle\nabla v,\nabla(uv')\rangle\Delta_pv
			+pw^{p/2-1}\langle\nabla v,\nabla(w^{p/2})\rangle\\
			&\quad-p(uv')w^{p-2}
			\bigl(|\nabla\nabla v|^2+\Ric(\nabla v,\nabla v)\bigr)
			-\frac{p(p-2)}4(uv')w^{p-3}|\nabla w|^2\\
			&\quad-\frac{p(p-2)}4(uv')w^{p-3}
			\left[|\nabla w|^2+
			\frac{p-2}w\langle\nabla v,\nabla w\rangle^2\right].
		\end{align*}
		Collecting the last two terms by \eqref{ms:eq:Hessian-A} gives
		\begin{equation}\label{ms:eq:evol-power-simplified}
			\begin{split}
				\square(w^{p/2})={}&pw^{p/2-1}\Delta_pv
				\langle\nabla v,\nabla(uv')\rangle
				+pw^{p/2-1}\langle\nabla v,\nabla(w^{p/2})\rangle\\
				&-p(uv')w^{p-2}
				\bigl(|\nabla\nabla v|_A^2+\Ric(\nabla v,\nabla v)\bigr).
			\end{split}
		\end{equation}
		Alternatively, \eqref{ms:eq:bochner-power} gives
		\eqref{ms:eq:evol-power} directly upon time differentiation of
		$w^{p/2}$; substituting
		$\nabla v_t=(uv')\nabla\Delta_pv+(\Delta_pv)\nabla(uv')
		+\nabla(w^{p/2})$
		then gives \eqref{ms:eq:evol-power-simplified}.
	\end{proof}

	Set $h=|\nabla\psi|^2$ and
	$\Box=\partial_t-\psi'\mathcal L_p$.
	Since $\psi=v$ and $\psi'=uv'$, these operators coincide with
	those in \eqref{ms:eq:Lp}--\eqref{ms:eq:square}.
	
	\begin{lemma}\label{ms:lem:psi-evolutions}
		On $\mathcal N$, we have
		\begin{align}
			\Box\psi_t
			&=ph^{p/2-1}\langle\nabla\psi,\nabla\psi_t\rangle
			+(\psi')_t\Delta_p\psi,
			\label{ms:eq:PF5}\\
			\Box h
			&=2\langle\nabla\psi,\nabla\psi_t\rangle
			-2\psi'h^{p/2-1}
			\Bigl(
			|\nabla\nabla\psi|^2+\Ric(\nabla\psi,\nabla\psi)
			+h^{1-p/2}\langle\nabla\psi,\nabla\Delta_p\psi\rangle
			\Bigr)\notag\\
			&-\left(\frac p2-1\right)\psi'h^{p/2-2}|\nabla h|^2,
			\label{ms:eq:PF6}\\
			\Box(h^{p/2})
			&=ph^{p/2-1}\langle\nabla\psi,\nabla\psi_t\rangle
			-p\psi'h^{p-2}
			\Bigl(
			|\nabla\nabla\psi|_A^2+\Ric(\nabla\psi,\nabla\psi)
			+h^{1-p/2}\langle\nabla\psi,\nabla\Delta_p\psi\rangle
			\Bigr).
			\label{ms:eq:PF7}
		\end{align}
	\end{lemma}
	
	\begin{proof}
		These identities follow directly from
		\eqref{ms:eq:evol-vt}--\eqref{ms:eq:evol-power}
		by setting $v=\psi$, $w=h$, and $uv'=\psi'$.
	\end{proof}

	\section{Li--Yau and Hamilton type estimates}
	\label{ms:grad:section}
	
	\subsection{Li--Yau type estimates}
	\label{ms:grad:LY}
	
	Throughout this subsection $(M^n,g)$ is closed, $\Ric\ge0$,
	$p\ge2$. We retain
	$\varphi=\log u$, $\psi(\varphi)=v(e^\varphi)$,
	$h=|\nabla\psi|^2$, and $b=\psi''/\psi'$.
	The multiplier in $\Box=\partial_t-\psi'\mathcal L_p$ is
	$\psi'=u^{(p-2)/(p-1)}F'(u)>0$.
	
	For $\alpha\ge1$, define the Harnack quantity
	\begin{equation}\label{ms:eq:H-def}
		H=h^{p/2}-\alpha\psi_t.
	\end{equation}
	It follows from \eqref{ms:eq:psi-equation} that
	\begin{equation}\label{ms:eq:H-laplacian}
		\Delta_p\psi
		=-\frac{H+(\alpha-1)h^{p/2}}{\alpha\psi'}.
	\end{equation}
	
	\begin{proposition}\label{ms:prop:H-evolution}
		On $\mathcal{N}$, we have
		\begin{align}
			\Box H={}&ph^{p/2-1}\langle\nabla\psi,\nabla H\rangle-p\psi'h^{p-2}
			\bigl(|\nabla\nabla\psi|_A^2+
			\Ric(\nabla\psi,\nabla\psi)\bigr)\\
			&-\frac{b}{\alpha\psi'}
			[H+(\alpha-1)h^{p/2}][H+(p-1)h^{p/2}].
			\label{ms:eq:H1}
		\end{align}
		
	\end{proposition}
	
	\begin{proof}
		Combining \eqref{ms:eq:PF5} and \eqref{ms:eq:PF7}, and using
		\eqref{ms:eq:H-laplacian}, we obtain
		\begin{align*}
			\Box H
			&=p(1-\alpha)h^{p/2-1}
			\langle\nabla\psi,\nabla\psi_t\rangle
			-p\psi'h^{p/2-1}
			\langle\nabla\psi,\nabla\Delta_p\psi\rangle-\alpha(\psi')_t\Delta_p\psi
			-p\psi'h^{p-2}
			\bigl(|\nabla\nabla\psi|_A^2
			+\Ric(\nabla\psi,\nabla\psi)\bigr)\\
			&=p(1-\alpha)h^{p/2-1}
			\langle\nabla\psi,\nabla\psi_t\rangle
			+\frac{p\psi'}{\alpha}h^{p/2-1}
			\left\langle\nabla\psi,
			\nabla\frac{H+(\alpha-1)h^{p/2}}{\psi'}\right\rangle\\
			&\quad+\frac{(\psi')_t}{\psi'}
			\bigl[H+(\alpha-1)h^{p/2}\bigr]
			-p\psi'h^{p-2}
			\bigl(|\nabla\nabla\psi|_A^2
			+\Ric(\nabla\psi,\nabla\psi)\bigr)\\
			&=ph^{p/2-1}\langle\nabla\psi,\nabla H\rangle
			-\frac{p}{\alpha\psi'}h^{p/2-1}
			\bigl[H+(\alpha-1)h^{p/2}\bigr]
			\langle\nabla\psi,\nabla\psi'\rangle\\
			&\quad+\frac{(\psi')_t}{\psi'}
			\bigl[H+(\alpha-1)h^{p/2}\bigr]
			-p\psi'h^{p-2}
			\bigl(|\nabla\nabla\psi|_A^2
			+\Ric(\nabla\psi,\nabla\psi)\bigr)\\
			&=ph^{p/2-1}\langle\nabla\psi,\nabla H\rangle
			-\frac{\psi''}{\alpha(\psi')^2}
			\bigl[H+(\alpha-1)h^{p/2}\bigr]
			\bigl[H+(p-1)h^{p/2}\bigr]\\
			&\quad-p\psi'h^{p-2}
			\bigl(|\nabla\nabla\psi|_A^2
			+\Ric(\nabla\psi,\nabla\psi)\bigr)\\
			&=ph^{p/2-1}\langle\nabla\psi,\nabla H\rangle
			-\frac{\psi''}{\alpha(\psi')^2}H^2
			-\frac{(p+\alpha-2)\psi''}{\alpha(\psi')^2}h^{p/2}H\\
			&\quad-\frac{(\alpha-1)(p-1)\psi''}
			{\alpha(\psi')^2}h^p
			-p\psi'h^{p-2}
			\bigl(|\nabla\nabla\psi|_A^2
			+\Ric(\nabla\psi,\nabla\psi)\bigr).
		\end{align*}
		
	\end{proof}
	
	\begin{proof}[\bf Proof of Theorem~\ref{ms:thm:LY}]
		Let $c\in C^1([\tau,T))$ satisfy the stated quadratic
		condition and $H(\cdot,\tau)\le c(\tau)$.
		By the Cauchy--Schwarz inequality applied to
		$A^{1/2}(\nabla\nabla\psi)A^{1/2}$,
		\begin{equation}\label{ms:eq:psi-trace}
			h^{p-2}|\nabla\nabla\psi|_A^2
			\ge\frac{(\Delta_p\psi)^2}{n}.
		\end{equation}
		Using $\Ric\ge0$ and
		\[
		H=-\alpha\psi'\Delta_p\psi-(\alpha-1)h^{p/2},
		\]
		we obtain from \eqref{ms:eq:H1}
		\begin{align}
			&\Box(H-c)
			-ph^{p/2-1}\langle\nabla\psi,\nabla(H-c)\rangle
			\notag\\
			&\quad\le
			-\frac{p\psi'}n(\Delta_p\psi)^2
			+\frac{\psi''}{\psi'}\Delta_p\psi
			\bigl[H+(p-1)h^{p/2}\bigr]-c'
			\notag\\
			&\quad=
			(p-\alpha)\psi''h^{p/2-1}
			\langle\nabla\psi,\nabla\varphi\rangle\Delta_p\psi
			-\left(\frac pn\psi'+\alpha\psi''\right)
			(\Delta_p\psi)^2-c'.
			\label{ms:eq:H-contact-inequality}
		\end{align}
		We examine a possible first contact point $(x_0,t_0)$
		where $H-c$ reaches zero from below. At this point,
		\[
		H=c(t_0),\qquad
		\nabla(H-c)=0,\qquad
		\partial_t(H-c)\ge0,\qquad
		\mathcal L_p(H-c)\le0.
		\]
		Consequently, \eqref{ms:eq:H-contact-inequality} gives
		\begin{align}
			0\le{}&
			(p-\alpha)\psi''h^{p/2-1}
			\langle\nabla\psi,\nabla\varphi\rangle\Delta_p\psi
			-\left(\frac pn\psi'+\alpha\psi''\right)
			(\Delta_p\psi)^2-c'(t_0).
			\label{ms:eq:H-contact}
		\end{align}
				The following pointwise calculation is carried out on
		$\{h>0\}$; for $p=2$ it holds everywhere.
		
		Since $\nabla\psi=\psi'\nabla\varphi$, we have
		\[
		h^{p/2}
		=\psi'h^{p/2-1}
		\langle\nabla\psi,\nabla\varphi\rangle.
		\]
		Moreover, at $(x_0,t_0)$,
		\begin{align*}
			\Delta_p\psi
			&=\frac{\psi_t-h^{p/2}}{\psi'}=\frac{(h^{p/2}-H)/\alpha-h^{p/2}}{\psi'}\\
			&=-\frac{H}{\alpha\psi'}
			+\frac{1-\alpha}{\alpha}
			h^{p/2-1}\langle\nabla\psi,\nabla\varphi\rangle=-\frac{c(t_0)}{\alpha\psi'}
			+\frac{1-\alpha}{\alpha}
			h^{p/2-1}\langle\nabla\psi,\nabla\varphi\rangle.
		\end{align*}
		Substituting into \eqref{ms:eq:H-contact}, we obtain
		\begin{align*}
			0\le{}&
			(p-\alpha)\psi''h^{p/2-1}
			\langle\nabla\psi,\nabla\varphi\rangle
			\left[
			-\frac{c}{\alpha\psi'}
			+\frac{1-\alpha}{\alpha}
			h^{p/2-1}\langle\nabla\psi,\nabla\varphi\rangle
			\right]\\
			&-\left(\frac pn\psi'+\alpha\psi''\right)
			\left[
			-\frac{c}{\alpha\psi'}
			+\frac{1-\alpha}{\alpha}
			h^{p/2-1}\langle\nabla\psi,\nabla\varphi\rangle
			\right]^2-c',
		\end{align*}
		where $c=c(t_0)$ and $c'=c'(t_0)$.
		
		Collecting the quadratic, linear, and constant terms yields
		\begin{equation}\label{ms:eq:H-contact-quadratic}
			0\le
			\mathcal A h^{p-2}
			\langle\nabla\psi,\nabla\varphi\rangle^2
			+\mathcal B h^{p/2-1}
			\langle\nabla\psi,\nabla\varphi\rangle
			+\mathcal C,
		\end{equation}
		where
		\begin{align*}
			\mathcal A
			&=\frac{(p-\alpha)(1-\alpha)}{\alpha}\psi''
			-\frac{(1-\alpha)^2}{\alpha^2}
			\left(\frac pn\psi'+\alpha\psi''\right)=\frac{1-\alpha}{\alpha^2}
			\left[
			\alpha(p-1)\psi''
			+\frac{p(\alpha-1)}n\psi'
			\right],\\
			\mathcal B
			&=-\frac{(p-\alpha)\psi''}{\alpha\psi'}c
			+\frac{2(1-\alpha)}{\alpha^2\psi'}
			\left(\frac pn\psi'+\alpha\psi''\right)c=-\left[
			\frac{2p(\alpha-1)}{n\alpha^2}
			+\frac{p+\alpha-2}{\alpha}\frac{\psi''}{\psi'}
			\right]c,\\
			\mathcal C
			&=-\left(\frac pn\psi'+\alpha\psi''\right)
			\frac{c^2}{\alpha^2(\psi')^2}-c'.
		\end{align*}
		
		Recalling that
		\[
		\psi'=u^{\frac{p-2}{p-1}}F'(u),\qquad
		\psi''=u^{\frac{p-2}{p-1}}
		\left[uF''(u)+\frac{p-2}{p-1}F'(u)\right],
		\]
		these coefficients can also be written as
		\begin{align*}
			\mathcal A
			&=\frac{1-\alpha}{\alpha^2}
			u^{\frac{p-2}{p-1}}
			\left\{
			\alpha(p-1)uF''(u)
			+\left[\alpha(p-2)+\frac{p(\alpha-1)}n\right]F'(u)
			\right\},\\
			\mathcal B
			&=-\left\{
			\frac{2p(\alpha-1)}{n\alpha^2}
			+\frac{p+\alpha-2}{\alpha}
			\left[\frac{uF''(u)}{F'(u)}+\frac{p-2}{p-1}\right]
			\right\}c,\\
			\mathcal C
			&=-\frac{
				\left\{\left[\frac pn+\frac{\alpha(p-2)}{p-1}\right]F'(u)
				+\alpha uF''(u)\right\}c^2
			}{
				\alpha^2u^{\frac{p-2}{p-1}}F'(u)^2
			}-c'.
		\end{align*}
		Since
		\[
		h^{p/2-1}\langle\nabla\psi,\nabla\varphi\rangle
		=\frac{h^{p/2}}{\psi'}\ge0,
		\]
		the quadratic condition \eqref{ms:eq:LY-polynomial}
		makes the right-hand side of
		\eqref{ms:eq:H-contact-quadratic} nonpositive.
		
		To complete the comparison also in the non-strict case,
		fix $\tau<t_1<T$. The same calculation with $H$ in place of $c$,
		together with the quadratic condition at $c$, gives
		$\Box(H-c)-ph^{p/2-1}\langle\nabla\psi,\nabla(H-c)\rangle
		\le C_0(H-c)$ where $H>c$, for some $C_0\ge0$ on
		$M\times[\tau,t_1]$, since the scalar reaction is a quadratic
		polynomial in $H$ with bounded coefficients on this strip.
		For $K>C_0$, a positive maximum of $e^{-K(t-\tau)}(H-c)$ occurs
		after $\tau$; for $p>2$, it lies in $\{h>0\}$ because $H=0$ at
		critical points of the smooth pressure, and for $p=2$, the operator
		is defined everywhere. At such a maximum,
		$0\le\Box\bigl(e^{-K(t-\tau)}(H-c)\bigr)
		\le(C_0-K)e^{-K(t-\tau)}(H-c)<0$,
		a contradiction. Hence $H\le c(t)$ on
		$M\times[\tau,t_1]$.
		Since $t_1<T$ is arbitrary,
		\[
		|\nabla\psi|^p-\alpha\psi_t\le c(t),
		\qquad \tau\le t<T.
		\]
	\end{proof}

\subsection{Aronson--B\'enilan type estimate}
	\label{ms:grad:AB}

When $b$ is nonincreasing, the scalar barrier in
Theorem~\ref{ms:thm:LY} can be replaced by a lower bound for
$\Delta_pv$ itself.

	We now work with the pressure $v$ and $E=\Delta_pv$. Since $v$ is
	strictly increasing in $u$, the elasticity $b$ in
	\eqref{ms:eq:b-def} is also a function of $v$; we write
	$b_v=\dd b/\dd v=b'(u)/v'(u)$. Along the solution,
	$\nabla(uv')=b\nabla v$ and $\nabla b=b_v\nabla v$, because
	$uv'=\psi'$ and $\nabla\psi'=(\psi''/\psi')\nabla\psi$.
	
	\begin{lemma}\label{ms:lem:AB-evolution}
		On $\mathcal N$, the function $E=\Delta_pv$ satisfies
		\begin{equation}\label{ms:eq:AB-evolution}
			\begin{aligned}
				\square E={}&\bigl[2(p-1)b+p\bigr]w^{p/2-1}
				\langle\nabla v,\nabla E\rangle+(p-1)bE^2
				+(p-1)b_vw^{p/2}E\\
				&+pw^{p-2}\bigl(|\nabla\nabla v|_A^2
				+\Ric(\nabla v,\nabla v)\bigr).
			\end{aligned}
		\end{equation}
	\end{lemma}
	
	\begin{proof}
		By \eqref{ms:eq:time-linearization} and
		\eqref{ms:eq:pressure-pde},
		$\partial_tE=\mathcal L_p\bigl((uv')E\bigr)+\mathcal L_p(w^{p/2})$.
		By \eqref{ms:eq:Lp-product} and $A\nabla v=(p-1)\nabla v$,
		\[
		\mathcal L_p\bigl((uv')E\bigr)=(uv')\mathcal L_pE
		+E\,\mathcal L_p(uv')+2(p-1)bw^{p/2-1}\langle\nabla v,\nabla E\rangle,
		\]
		and
		\[
		\mathcal L_p(uv')=\Div\bigl(w^{p/2-1}A(b\nabla v)\bigr)
		=(p-1)\Div\bigl(bw^{p/2-1}\nabla v\bigr)=(p-1)\bigl(bE+b_vw^{p/2}\bigr).
		\]
		Adding \eqref{ms:eq:bochner-power} and moving
		$(uv')\mathcal L_pE$ to the left-hand side gives
		\eqref{ms:eq:AB-evolution}.
	\end{proof}
	
	\begin{proof}[\bf Proof of Theorem~\ref{ms:thm:AB}]
		By \eqref{ms:eq:psi-trace}, written for $v$, and $\Ric\ge0$,
		Lemma~\ref{ms:lem:AB-evolution} gives on $\mathcal N$
		\begin{equation}\label{ms:eq:AB-inequality}
			\square E\ge\bigl[2(p-1)b+p\bigr]w^{p/2-1}
			\langle\nabla v,\nabla E\rangle+\varrho E^2
			+(p-1)b_vw^{p/2}E,
		\end{equation}
		where we used $(p-1)b+p/n\ge\varrho$. Fix $\varepsilon\in(0,T)$
		and $T'\in(\varepsilon,T)$, and set $Z=(t-\varepsilon)E$ on
		$M\times[\varepsilon,T']$, so that $Z(\cdot,\varepsilon)=0$.
		Suppose that $\min Z<-1/\varrho$. For $\eta>0$ small, let $t_1$ be
		the first time at which $\min_MZ(\cdot,t_1)=-(1+\eta)/\varrho$,
		attained at $x_1$. For $p>2$, $E=0$ at critical points of the
		smooth pressure, since $\nabla\bigl(w^{p/2-1}\nabla v\bigr)=w^{p/2-1}A(\nabla\nabla v)$
		vanishes there; for $p=2$ the operator is defined everywhere.
		Hence $(x_1,t_1)$ is a point of $\mathcal N$ when $p>2$, and
		\[
		\nabla E=0,\qquad
		\mathcal L_pE=w^{p/2-1}\operatorname{tr}\bigl(A\nabla\nabla E\bigr)\ge0,
		\qquad\partial_tZ\le0.
		\]
		Since $b_v\le0$ and $E<0$, the last term in
		\eqref{ms:eq:AB-inequality} is nonnegative. Therefore
		\[
		0\ge\partial_tZ=E+(t_1-\varepsilon)\partial_tE
		\ge E+(t_1-\varepsilon)\varrho E^2
		=\frac{Z}{t_1-\varepsilon}(1+\varrho Z)
		=\frac{(1+\eta)\eta}{\varrho(t_1-\varepsilon)}>0,
		\]
		a contradiction. Thus $(t-\varepsilon)\Delta_pv\ge-1/\varrho$ on
		$M\times[\varepsilon,T']$, and letting $\varepsilon\to0$,
		$T'\to T$ proves the first inequality in \eqref{ms:eq:AB-bound}.
		The second follows from $\psi=v$, $\psi'=uv'$ and
		\eqref{ms:eq:pressure-pde}, which give
		$|\nabla\psi|^p-\psi_t=-\psi'\Delta_pv$.
	\end{proof}
	
	Integrating \eqref{ms:eq:LY-barrier} or \eqref{ms:eq:AB-bound} along
	minimizing geodesics gives Harnack inequalities for $u$ in the usual
	way \cite{LY}; we omit the details.

	\subsection{Hamilton type estimate}
	
	In this subsection, $p>1$ and $h=|\nabla\psi|^2$. The calculation takes
	place on $\mathcal N=\{h>0\}$.
	
	\begin{lemma}[\cite{Xu}]\label{ms:lem:Xu}
		Let $S(n)$ denote the vector space of real symmetric
		$n\times n$ matrices, and let $|H|$ denote the Frobenius norm.
		For any $\alpha,\beta\in\mathbb R$,
		\begin{equation}\label{ms:eq:Xu-matrix}
			\max_{\substack{H\in S(n)\setminus\{0\}\\ |e|=1}}
			\left[\left(\frac{\alpha H+\beta(\tr H)I_n}{|H|}\right)(e,e)\right]^2
			=(\alpha+\beta)^2+(n-1)\beta^2,
		\end{equation}
		where $I_n$ is the identity matrix.
	\end{lemma}

	Define another Harnack quantity by
	\begin{equation}\label{ms:eq:G-def}
		G:=\frac{|\nabla\psi|^p}{(\alpha-\psi)^q},
	\end{equation}
	where $q>0$ is a constant to be chosen later.
	
	\begin{proposition}\label{ms:prop:BG}
		On $\mathcal N$, we have
		\begin{align}
			\Box G
			&=\bigl[pb-q(p-2)d\bigr]G\Delta_p\psi
			+q(p-1)\bigl[(q-1)d-1\bigr](\alpha-\psi)^{q-1}G^2
			\notag\\
			&\quad+\bigl[p-2q(p-1)d\bigr]|\nabla\psi|^{p-2}
			\langle\nabla G,\nabla\psi\rangle
			-\frac{p\psi'h^{p-2}}{(\alpha-\psi)^q}
			\bigl(|\nabla\nabla\psi|_A^2
			+\Ric(\nabla\psi,\nabla\psi)\bigr),
			\label{ms:eq:BG}
		\end{align}
		where
		\[
		b=\frac{\psi''}{\psi'},\qquad
		d=\frac{\psi'}{\alpha-\psi}.
		\]
	\end{proposition}
	
	\begin{proof}
		For smooth functions $f$ and $g>0$, the quotient rule gives
		\begin{equation}\label{ms:eq:quotient}
			\Box\left(\frac fg\right)
			=\frac1g\Box f-\frac f{g^2}\Box g
			+2\psi'h^{p/2-1}
			\left\langle
			A\nabla\left(\frac fg\right),\nabla\log g
			\right\rangle.
		\end{equation}
		Applying this identity to \eqref{ms:eq:G-def}, we obtain
		\begin{align}
			\Box G
			&=\frac1{(\alpha-\psi)^q}\Box|\nabla\psi|^p
			-\frac{|\nabla\psi|^p}{(\alpha-\psi)^{2q}}
			\Box(\alpha-\psi)^q
			+2\psi'|\nabla\psi|^{p-2}
			\left\langle
			A\nabla G,\nabla\log\bigl((\alpha-\psi)^q\bigr)
			\right\rangle.
			\label{ms:eq:boxdot1}
		\end{align}
		Note that
		\[
		\nabla\psi'=\frac{\psi''}{\psi'}\nabla\psi,
		\qquad
		\nabla|\nabla\psi|^p
		=(\alpha-\psi)^q\nabla G
		-q(\alpha-\psi)^{q-1}G\nabla\psi.
		\]
		Using \eqref{ms:eq:PF7} and the pressure equation
		\eqref{ms:eq:psi-equation}, we have
		\begin{align}
			\frac1{(\alpha-\psi)^q}\Box|\nabla\psi|^p
			&=\frac1{(\alpha-\psi)^q}
			\left[
			p\frac{\psi''}{\psi'}|\nabla\psi|^p\Delta_p\psi
			+p(\alpha-\psi)^q|\nabla\psi|^{p-2}
			\langle\nabla G,\nabla\psi\rangle
			\right]
			\notag\\
			&\quad-\frac1{(\alpha-\psi)^q}
			\left[
			pq|\nabla\psi|^p(\alpha-\psi)^{q-1}G
			+p\psi'h^{p-2}
			\bigl(|\nabla\nabla\psi|_A^2
			+\Ric(\nabla\psi,\nabla\psi)\bigr)
			\right]\label{ms:eq:boxdot2}
			\\
			&=p\frac{\psi''}{\psi'}G\Delta_p\psi
			+p|\nabla\psi|^{p-2}
			\langle\nabla G,\nabla\psi\rangle
			-pq(\alpha-\psi)^{q-1}G^2
			-\frac{p\psi'h^{p-2}}{(\alpha-\psi)^q}
			\bigl(|\nabla\nabla\psi|_A^2
			+\Ric(\nabla\psi,\nabla\psi)\bigr).
			\notag
		\end{align}
		On the other hand, using
		$A\nabla\psi=(p-1)\nabla\psi$, we obtain
		\begin{align}
			\frac{|\nabla\psi|^p}{(\alpha-\psi)^{2q}}
			\Box(\alpha-\psi)^q
			&=
			\frac{|\nabla\psi|^p}{(\alpha-\psi)^{2q}}
			\left[
			\partial_t(\alpha-\psi)^q
			-\psi'\mathcal L_p(\alpha-\psi)^q
			\right]
			\notag\\
			&\quad=
			\frac{|\nabla\psi|^p}{(\alpha-\psi)^{2q}}
			\left[
			-q(\alpha-\psi)^{q-1}
			\bigl(\psi'\Delta_p\psi+|\nabla\psi|^p\bigr)
			+q(p-1)\psi'(\alpha-\psi)^{q-1}\Delta_p\psi
			\right]
			\notag\\
			&\qquad-
			\frac{|\nabla\psi|^p}{(\alpha-\psi)^{2q}}
			q(q-1)\psi'(\alpha-\psi)^{q-2}|\nabla\psi|^{p-2}
			\langle\nabla\psi,A\nabla\psi\rangle
			\notag\\
			&\quad=
			q(p-2)\frac{\psi'}{\alpha-\psi}G\Delta_p\psi
			-q\left[
			1+(p-1)(q-1)\frac{\psi'}{\alpha-\psi}
			\right](\alpha-\psi)^{q-1}G^2.
			\label{ms:eq:boxdot3}
		\end{align}		
		For the remaining term, we have
		\begin{align}
			2\psi'|\nabla\psi|^{p-2}
			\left\langle
			A\nabla G,\nabla\log\bigl((\alpha-\psi)^q\bigr)
			\right\rangle
			&=
			-\frac{2q\psi'}{\alpha-\psi}|\nabla\psi|^{p-2}
			\langle A\nabla G,\nabla\psi\rangle
			\notag\\
			&\quad=
			-2q(p-1)\frac{\psi'}{\alpha-\psi}|\nabla\psi|^{p-2}
			\langle\nabla G,\nabla\psi\rangle.
			\label{ms:eq:boxdot4}
		\end{align}
		Substituting \eqref{ms:eq:boxdot2}--\eqref{ms:eq:boxdot4}
		into \eqref{ms:eq:boxdot1}, and observing that
		\[
		-pq+q\bigl[1+(p-1)(q-1)d\bigr]
		=q(p-1)\bigl[(q-1)d-1\bigr],
		\]
		we obtain \eqref{ms:eq:BG}.
	\end{proof}
	
	\begin{lemma}\label{ms:prop:G-optimized}
		Let $p>1$ and suppose that $\Ric\ge-\kappa g$,
		where $\kappa\ge0$. Then, on $\mathcal N$,
		\begin{equation}\label{ms:eq:G-optimized}
			\begin{aligned}
				\Box G
				\le{}&
				\Biggl\{
				\frac{n-1}{4pd}[pb-q(p-2)d]^2
				-q(p-1)\left[
				1+b+\left(1-\frac{q(p-1)}p\right)d
				\right]
				\Biggr\}(\alpha-\psi)^{q-1}G^2\\
				&+\bigl[p+(p-1)b-q(p-1)d\bigr]
				h^{p/2-1}\langle\nabla G,\nabla\psi\rangle
				+p\kappa\psi'h^{p/2-1}G.
			\end{aligned}
		\end{equation}
	\end{lemma}
	
	\begin{proof}
		By the definition of $G$,
		\begin{equation}\label{ms:eq:G-directional}
			h^{p/2-1}\langle\nabla G,\nabla\psi\rangle
			=
			\frac{ph^{p-2}}{(\alpha-\psi)^q}
			\nabla\nabla\psi(\nabla\psi,\nabla\psi)
			+q(\alpha-\psi)^{q-1}G^2.
		\end{equation}
		We introduce a scalar parameter $\eta$
		and use this identity to rearrange the first-order term.
		Combining \eqref{ms:eq:BG} with
		\eqref{ms:eq:G-directional}, we obtain
		\begin{equation}\label{ms:eq:G-Xu-parameter}
			\begin{aligned}
				\Box G
				={}&-\frac{p\psi'h^{p-2}}{(\alpha-\psi)^q}
				|\nabla\nabla\psi|_A^2+G\left\{
				[pb-q(p-2)d]\Delta_p\psi
				-p\eta h^{p/2-2}
				\nabla\nabla\psi(\nabla\psi,\nabla\psi)
				\right\}\\
				&+\left\{q(p-1)[(q-1)d-1]-q\eta\right\}
				(\alpha-\psi)^{q-1}G^2+\bigl[p-2q(p-1)d+\eta\bigr]
				h^{p/2-1}\langle\nabla G,\nabla\psi\rangle\\
				&-\frac{p\psi'h^{p-2}}{(\alpha-\psi)^q}
				\Ric(\nabla\psi,\nabla\psi).
			\end{aligned}
		\end{equation}
		
		Applying Lemma~\ref{ms:lem:Xu} to the symmetric matrix
		$A^{1/2}(\nabla\nabla\psi)A^{1/2}$ gives
		\begin{equation}\label{ms:eq:Xu-Hessian-bound}
			\begin{aligned}
				&-\frac{p\psi'h^{p-2}}{(\alpha-\psi)^q}
				|\nabla\nabla\psi|_A^2+
				G\left\{
				[pb-q(p-2)d]\Delta_p\psi
				-p\eta h^{p/2-2}
				\nabla\nabla\psi(\nabla\psi,\nabla\psi)
				\right\}\\
				&\quad\le
				\frac{(\alpha-\psi)^{q-1}G^2}{4pd}
				\left\{
				\left[pb-q(p-2)d-\frac{p\eta}{p-1}\right]^2
				+(n-1)[pb-q(p-2)d]^2
				\right\}.
			\end{aligned}
		\end{equation}
		The terms depending on $\eta$ satisfy
		\[
		\begin{aligned}
			&\frac1{4pd}
			\left[pb-q(p-2)d-\frac{p\eta}{p-1}\right]^2-q\eta=
			\frac{p}{4d(p-1)^2}
			\bigl[\eta-(p-1)(b+qd)\bigr]^2
			-q(p-1)b-\frac{q^2(p-1)}p d.
		\end{aligned}
		\]
		We therefore choose
		\[
		\eta=(p-1)(b+qd).
		\]
		Substituting this choice into
		\eqref{ms:eq:G-Xu-parameter} and
		\eqref{ms:eq:Xu-Hessian-bound}, and using
		\[
		-\frac{p\psi'h^{p-2}}{(\alpha-\psi)^q}
		\Ric(\nabla\psi,\nabla\psi)
		\le p\kappa\psi'h^{p/2-1}G,
		\]
		proves \eqref{ms:eq:G-optimized}.
	\end{proof}

	\begin{proof}[\bf Proof of Theorem~\ref{ms:thm:Hamilton}]
		We first bound the pressure gap independently of $\alpha$.
		Choose $L\ge\delta$ such that
		\[
		\frac{q^2}{pL}
		\left[(p-1)^2+\frac{(n-1)(p-2)^2}{4}\right]
		\max_{\mathcal U}\psi'\le\frac\beta2,
		\qquad
		\widetilde\alpha=\min\{\alpha,\sup_{\mathcal U}v+L\}.
		\]
		Then $\delta\le\widetilde\alpha-v\le L+\operatorname{osc}_{\mathcal U}v$.
		If $\widetilde\alpha<\alpha$, replacing $\alpha$ by
		$\widetilde\alpha$ increases $d$, and the new $d$ is at most
		$\max_{\mathcal U}\psi'/L$.
		In the expanded structural condition
		\eqref{ms:eq:Hamilton-condition}, the negative reciprocal term
		increases, while the linear term in $d$ decreases by at most
		$\beta/2$, by the choice of $L$.
		Thus the modified structural expression is at least $\beta/2$.
		The corresponding quotient dominates $G$, so it suffices to
		prove the estimate for this modified upper level.
		We retain the notation $\alpha$ and $\beta$ for the modified
		level and positive structural lower bound.
		All coefficient and pressure-gap bounds now depend only on
		$n,p,q,\delta,\beta$ and $\Lambda$: indeed
		$\delta\le\alpha-\psi\le L+\Lambda$, $0<d\le\Lambda/\delta$, and
		$\psi'$, $|b|$ are bounded by $\Lambda$.
		
		Let $\Psi=\Psi(x,t)$ be a space--time cutoff supported
		spatially in $B(x_0,R)$ and vanishing near the initial
		time of $Q_{R,T}$, such that
		\[
		0\le\Psi\le1,
		\qquad
		\Psi=1\quad\text{on }Q_{R/2,T/2}.
		\]
		Using a sufficiently high power of a radial cutoff
		and a temporal cutoff, we may arrange that
		\begin{equation}\label{ms:eq:Hamilton-cutoff}
			|\nabla\Psi|
			\le\frac{C}{R}\Psi^{1-1/p},
			\qquad
			|\Psi_t|\le\frac{C}{T}.
		\end{equation}
		Under the sectional curvature bound
		$\operatorname{Sec}\ge-\kappa$, Hessian comparison also gives
		\begin{equation}\label{ms:eq:Hamilton-cutoff-Hessian}
			-\operatorname{tr}(A\nabla\nabla\Psi)
			\le
			C\left(\frac1{R^2}+\frac{\sqrt\kappa}{R}\right)
			\Psi^{1-2/p}.
		\end{equation}
		As in \cite{KN}, the calculation at the cut locus
		is justified by the supporting cutoff argument.
		
		Suppose that $\Psi G$ attains a positive maximum at
		$(x_1,t_1)$. Then $\Psi>0$, $h>0$, and
		\begin{equation}\label{ms:eq:Hamilton-critical-relations}
			\nabla G=-\frac{G}{\Psi}\nabla\Psi,
			\qquad
			\nabla h
			=-\frac{2h}{p}
			\left(
			\frac{\nabla\Psi}{\Psi}
			+\frac{q\nabla\psi}{\alpha-\psi}
			\right).
		\end{equation}
		At this point, the maximum principle gives
		\begin{equation}\label{ms:eq:Hamilton-product-maximum}
			\begin{aligned}
				0\le\Box(\Psi G)
				&=\Psi\Box G+G\Box\Psi
				-2\psi'h^{p/2-1}
				\langle A\nabla\Psi,\nabla G\rangle\\
				&=\Psi\Box G+G\Box\Psi
				+2\psi'h^{p/2-1}\frac{G}{\Psi}
				\langle A\nabla\Psi,\nabla\Psi\rangle.
			\end{aligned}
		\end{equation}
		
		We first compute the terms arising from the divergence
		form of the operator. We have
		\begin{equation}\label{ms:eq:Hamilton-cutoff-operator}
			\begin{aligned}
				\Box\Psi
				={}&\Psi_t
				-\psi'h^{p/2-1}
				\operatorname{tr}(A\nabla\nabla\Psi)-\psi'
				\left\langle
				\operatorname{div}(h^{p/2-1}A),\nabla\Psi
				\right\rangle.
			\end{aligned}
		\end{equation}
		Direct differentiation yields
		\begin{align*}
			\left\langle
			\operatorname{div}(h^{p/2-1}A),\nabla\Psi
			\right\rangle=
			(p-2)h^{p/2-2}
			\left[
			\langle\nabla h,\nabla\Psi\rangle
			-\frac{
				\langle\nabla h,\nabla\psi\rangle
				\langle\nabla\Psi,\nabla\psi\rangle
			}{h}
			\right]+
			(p-2)\frac{\Delta_p\psi}{h}
			\langle\nabla\Psi,\nabla\psi\rangle.
		\end{align*}
		Using \eqref{ms:eq:Hamilton-critical-relations},
		this becomes
		\begin{equation}\label{ms:eq:Hamilton-cutoff-divergence}
			\begin{aligned}
				\left\langle
				\operatorname{div}(h^{p/2-1}A),\nabla\Psi
				\right\rangle=
				-\frac{2(p-2)}{p\Psi}h^{p/2-1}
				\left[
				|\nabla\Psi|^2
				-\frac{\langle\nabla\Psi,\nabla\psi\rangle^2}{h}
				\right]+
				(p-2)\frac{\Delta_p\psi}{h}
				\langle\nabla\Psi,\nabla\psi\rangle.
			\end{aligned}
		\end{equation}
		
		We combine the last term in
		\eqref{ms:eq:Hamilton-cutoff-divergence}
		with the Hessian terms in
		\eqref{ms:eq:G-Xu-parameter}, taking
		$\eta=(p-1)(b+qd)$.
		The resulting Hessian contribution to
		\eqref{ms:eq:Hamilton-product-maximum} is
		\begin{equation}\label{ms:eq:Hamilton-local-Hessian}
			\begin{aligned}
				-\frac{p\psi'\Psi h^{p-2}}{(\alpha-\psi)^q}
				|\nabla\nabla\psi|_A^2+
				\Psi G\Biggl\{
				\left[
				pb-q(p-2)d
				-\frac{(p-2)\psi'}{h\Psi}
				\langle\nabla\Psi,\nabla\psi\rangle
				\right]\Delta_p\psi\Biggr.\\
				\Biggl.\qquad\qquad
				-p(p-1)(b+qd)h^{p/2-2}
				\nabla\nabla\psi(\nabla\psi,\nabla\psi)
				\Biggr\}.
			\end{aligned}
		\end{equation}
		Applying the same matrix estimate as in
		\eqref{ms:eq:Xu-Hessian-bound} and expanding the
		resulting squares, we bound
		\eqref{ms:eq:Hamilton-local-Hessian} from above by
		\begin{equation}\label{ms:eq:Hamilton-local-matrix-bound}
			\begin{aligned}
				&\frac{\Psi(\alpha-\psi)^{q-1}G^2}{4pd}
				\left\{
				4q^2(p-1)^2d^2
				+(n-1)[pb-q(p-2)d]^2
				\right\}\\
				&\quad-
				\frac{p-2}{2p}
				\left\{
				p(n-1)b-q[n(p-2)+p]d
				\right\}
				h^{p/2-1}G\langle\nabla\Psi,\nabla\psi\rangle+
				\frac{n(p-2)^2\psi'}{4p\Psi}
				h^{p/2-2}G
				\langle\nabla\Psi,\nabla\psi\rangle^2.
			\end{aligned}
		\end{equation}
		The first line is the same quadratic contribution as
		in the preceding lemma, multiplied by $\Psi$.
		The remaining two lines are bounded above by
		\[
		C h^{(p-1)/2}G|\nabla\Psi|
		+C\psi'h^{p/2-1}G\frac{|\nabla\Psi|^2}{\Psi}.
		\]
		
		The sectional curvature bound gives
		$\Ric\ge-(n-1)\kappa g$.
		By the structural assumption, the reaction coefficient is at most
		$-\beta(\alpha-\psi)^{q-1}$.
		Replacing $\beta$ by
		$\beta\min_{\mathcal U}(\alpha-v)^{q-1}>0$
		and retaining its notation, we obtain a uniform upper bound
		$-\beta$ for the coefficient of $G^2$.
		Combining
		\eqref{ms:eq:Hamilton-product-maximum}--%
		\eqref{ms:eq:Hamilton-local-matrix-bound},
		we therefore obtain at $(x_1,t_1)$
		\begin{equation}\label{ms:eq:Hamilton-local-maximum}
			\begin{aligned}
				\beta\Psi G^2
				\le{}&
				G|\Psi_t|
				-\psi'h^{p/2-1}G
				\operatorname{tr}(A\nabla\nabla\Psi)
				+p(n-1)\kappa\psi'\Psi h^{p/2-1}G\\
				&+C h^{(p-1)/2}G|\nabla\Psi|
				+C\psi'h^{p/2-1}G
				\frac{|\nabla\Psi|^2}{\Psi}.
			\end{aligned}
		\end{equation}
		Here the first-order term from
		\eqref{ms:eq:G-Xu-parameter} has been estimated using
		\[
		\begin{aligned}
			&\Psi\bigl[p+(p-1)b-q(p-1)d\bigr]
			h^{p/2-1}\langle\nabla G,\nabla\psi\rangle\\
			&\quad=
			-\bigl[p+(p-1)b-q(p-1)d\bigr]
			h^{p/2-1}G\langle\nabla\Psi,\nabla\psi\rangle,
		\end{aligned}
		\]
		and all terms quadratic in $\nabla\Psi$ have been
		included in the last term of
		\eqref{ms:eq:Hamilton-local-maximum}.
		
		Since
		\[
		h^{p/2-1}G
		=(\alpha-\psi)^{q(p-2)/p}G^{2-2/p},
		\qquad
		h^{(p-1)/2}G
		=(\alpha-\psi)^{q(p-1)/p}G^{2-1/p},
		\]
		we multiply \eqref{ms:eq:Hamilton-local-maximum}
		by $\Psi$ and apply
		\eqref{ms:eq:Hamilton-cutoff} and
		\eqref{ms:eq:Hamilton-cutoff-Hessian}.
		The fixed coefficient bounds on $\mathcal U$ then give
		\begin{equation}\label{ms:eq:Hamilton-cutoff-powers}
			\begin{aligned}
				\beta(\Psi G)^2
				\le C\Biggl[
				\frac{\Psi G}{T}
				+\frac{(\Psi G)^{2-1/p}}{R}
				+\left(
				\kappa+\frac1{R^2}+\frac{\sqrt\kappa}{R}
				\right)(\Psi G)^{2-2/p}
				\Biggr].
			\end{aligned}
		\end{equation}
		
		We estimate these three terms as in Xu's localization
		argument. Since $p>1$, Young's inequality yields
		\begin{align*}
			\frac{C}{T}\Psi G
			&\le
			\frac{\beta}{6}(\Psi G)^2+\frac{C}{T^2},\\
			\frac{C}{R}(\Psi G)^{2-1/p}
			&\le
			\frac{\beta}{6}(\Psi G)^2+\frac{C}{R^{2p}},\\
			C\left(
			\kappa+\frac1{R^2}+\frac{\sqrt\kappa}{R}
			\right)(\Psi G)^{2-2/p}
			&\le
			\frac{\beta}{6}(\Psi G)^2
			+C\left(
			\kappa+\frac1{R^2}+\frac{\sqrt\kappa}{R}
			\right)^p.
		\end{align*}
		Substituting into
		\eqref{ms:eq:Hamilton-cutoff-powers}, we conclude that
		\[
		(\Psi G)(x_1,t_1)
		\le
		C\left[
		\frac1T+\frac1{R^p}
		+\left(\frac{\sqrt\kappa}{R}\right)^{p/2}
		+\kappa^{p/2}
		\right].
		\]
		Since $(x_1,t_1)$ is a maximum point and $\Psi=1$
		on $Q_{R/2,T/2}$, the desired local estimate follows.
		By the initial reduction, the constant depends only on
		$n,p,q,\delta,\beta$ and $\Lambda$.
		
		For $p=2$, the divergence term in
		\eqref{ms:eq:Hamilton-cutoff-divergence} vanishes,
		and the spatial cutoff enters only through
		$\Delta\Psi$.
		Laplacian comparison under $\Ric\ge-\kappa g$
		therefore suffices, and the same argument gives
		\[
		G\le C\left(\frac1T+\frac1{R^2}+\kappa\right)
		\quad\text{on }Q_{R/2,T/2}.
		\]
	\end{proof}
	
	On closed manifolds the same maximum-principle argument
gives a time-decay estimate. On complete manifolds the local estimate
also yields a Liouville conclusion for ancient solutions.

\begin{corollary}[Global estimate and Liouville theorem]
		\label{ms:cor:Hamilton}
		Let $p>1$, and let $F$, $\mathcal U$, $\psi$, $q$ and $\beta$
		satisfy the structural hypotheses of
		Theorem~\ref{ms:thm:Hamilton} for some $\alpha>\sup_{\mathcal U}v$.
		\begin{enumerate}[label=\rm(\roman*)]
			\item If $(M^n,g)$ is closed with $\Ric\ge0$ and $u$ solves
			\eqref{ms:intro:PF} on $M\times(0,T)$ with values in
			$\mathcal U$, then
			\begin{equation}\label{ms:eq:Hamilton-global}
				t|\nabla\psi|^p\le
				\frac{(\alpha-\psi)^q}
				{\beta\min_{\mathcal U}(\alpha-v)^{q-1}}
				\qquad\text{on }M\times(0,T).
			\end{equation}
			\item If $(M^n,g)$ is connected and complete with
			$\operatorname{Sec}\ge0$ (with $\Ric\ge0$ when $p=2$), then every
			ancient solution of \eqref{ms:intro:PF} on $M\times(-\infty,0]$
			with values in $\mathcal U$ is constant.
		\end{enumerate}
	\end{corollary}

\begin{proof}
		(i) Fix $\varepsilon\in(0,T)$ and $T'\in(\varepsilon,T)$, and
		suppose that $(t-\varepsilon)G$ attains a positive maximum on
		$M\times[\varepsilon,T']$ at $(x_1,t_1)$. Then $t_1>\varepsilon$,
		$h>0$, $\nabla G=0$,
		$\mathcal L_pG=h^{p/2-1}\operatorname{tr}(A\nabla\nabla G)\le0$ and
		$\partial_t\bigl((t-\varepsilon)G\bigr)\ge0$ there. Since the
		structural expression is at least $\beta$, Lemma~\ref{ms:prop:G-optimized}
		with $\kappa=0$ gives
		$\Box G\le-\beta(\alpha-\psi)^{q-1}G^2$ at this point. Hence
		\[
		0\le\Box\bigl((t-\varepsilon)G\bigr)
		\le G-\beta(t_1-\varepsilon)(\alpha-\psi)^{q-1}G^2,
		\]
		so $(t_1-\varepsilon)G\le1/(\beta\min_{\mathcal U}(\alpha-v)^{q-1})$.
		Letting $\varepsilon\to0$ and $T'\to T$ gives
		\eqref{ms:eq:Hamilton-global}.
		
		(ii) Let $(x_0,t_0)\in M\times(-\infty,0]$. For every $R,T>0$, the
		solution is defined on $Q_{R,T}$, and Theorem~\ref{ms:thm:Hamilton}
		with $\kappa=0$ gives
		\[
		|\nabla\psi|(x_0,t_0)\le C\left(\frac1R+\frac1{T^{1/p}}\right)
		\bigl(\alpha-\psi(x_0,t_0)\bigr)^{q/p},
		\]
		with $C$ independent of $R$ and $T$. Letting $R,T\to\infty$ shows
		$\nabla\psi\equiv0$. Thus $u$ is spatially constant, and
		\eqref{ms:intro:PF} gives $u_t=0$.
	\end{proof}

	\section{\texorpdfstring{$\mathcal W$}{W}-entropy and entropy-power concavity}
	\label{ms:ent:section}
	\label{ms:rig:section}
	
	Throughout this section, $(M^n,g)$ is closed with $\Ric\ge0$,
	and $u>0$ is a solution of \eqref{ms:intro:PF}.
	We use the hypotheses and normalization of
	Theorem~\ref{ms:thm:entropy}.
	
	\subsection{Entropy variations and the \texorpdfstring{$\mathcal W$}{W}-formula}
	\label{ms:ent:variables}
	
	Recall that
	\begin{equation}\label{ms:ent:Upp}
		U'(s)=v(s),\qquad
		U''(s)=s^{-1/(p-1)}F'(s),\qquad
		U(0)=v(0)=0,
	\end{equation}
	and
	\[
	P(s)=sU'(s)-U(s).
	\]
	In particular,
	\begin{equation}\label{ms:ent:pressure-identities}
		P'(u)=uU''(u),\qquad
		\nabla P(u)=u\nabla v.
	\end{equation}
	
	Retaining $w=|\nabla v|^2$ and the matrix $A$ from the
	pressure calculations, we have
\begin{equation}\label{ms:ent:flux}
	\Delta_pv=\Div\bigl(w^{p/2-1}\nabla v\bigr).
\end{equation}
	The density and pressure equations become
\begin{equation}\label{ms:ent:equations}
	\begin{gathered}
		u_t=\Div\bigl(uw^{p/2-1}\nabla v\bigr),\qquad
		v_t=P'(u)\Delta_pv+w^{p/2},\\
		u_t=u\Delta_pv+\frac{u}{P'(u)}w^{p/2}.
	\end{gathered}
\end{equation}
	On $\{w>0\}$,
\begin{equation}\label{ms:ent:critical-square}
	\nabla\bigl(w^{p/2-1}\nabla v\bigr)
	=w^{p/2-1}A(\nabla\nabla v).
\end{equation}
	We use the weighted Hilbert--Schmidt norm
	$|\nabla\nabla v|_A$ defined in \eqref{ms:eq:Hessian-A}.
	At critical points, all weighted expressions below are
	interpreted by continuous extension through the pressure flux
	$w^{p/2-1}\nabla v$ and its derivative.
	
	Recall the internal energy and its dissipation:
	\begin{equation}\label{ms:ent:EI}
		\mathcal E(t)=\int_M U(u)\,\dd V,\qquad
		\mathcal I(t)=\int_M u|\nabla v|^p\,\dd V.
	\end{equation}
	Inside integrals, $U$ and $P$ are evaluated at $u$.
	
	\begin{lemma}[Entropy variations]
		\label{ms:ent:first-variation}
		\label{ms:ent:Bochner}
		The energy and its dissipation satisfy
		\begin{equation}\label{ms:ent:Eprime}
			\mathcal E'=-\mathcal I,\qquad
			\int_M P\Delta_pv\,\dd V=-\mathcal I,
		\end{equation}
		and
		\begin{equation}\label{ms:ent:Iprime}
			\mathcal I'
			=-p\int_M
			\left\{
			P\left[
			w^{p-2}\bigl(|\nabla\nabla v|_A^2
			+\Ric(\nabla v,\nabla v)\bigr)
			\right]
			+\bigl(uP'(u)-P\bigr)(\Delta_pv)^2
			\right\}\,\dd V.
		\end{equation}
	\end{lemma}
	
	\begin{proof}
		Integration by parts gives
		\[
		\mathcal E'
		=\int_M v\Div\bigl(uw^{p/2-1}\nabla v\bigr)\,\dd V
		=-\int_Muw^{p/2}\,\dd V
		=-\mathcal I.
		\]
		Similarly, $\nabla P=u\nabla v$ implies
		\[
		\int_M P\Delta_pv\,\dd V
		=-\int_Mw^{p/2-1}\langle\nabla P,\nabla v\rangle\,\dd V
		=-\mathcal I.
		\]
		
		Differentiating $\mathcal I$ and using
		\eqref{ms:ent:equations}, we obtain
		\begin{align}
			\mathcal I'
			&=
			\int_Mu_tw^{p/2}\,\dd V
			+p\int_Muw^{p/2-1}\langle\nabla v,\nabla v_t\rangle\,\dd V
			\notag\\
			&=
			\int_Mu_tw^{p/2}\,\dd V
			-p\int_Mu_tv_t\,\dd V
			\notag\\
			&=
			-p\int_M P'(u)u_t\Delta_pv\,\dd V
			-(p-1)\int_Mw^{p/2}u_t\,\dd V
			\notag\\
			&=
			-p\int_MuP'(u)(\Delta_pv)^2\,\dd V
			-p\int_Muw^{p/2}\Delta_pv\,\dd V
			-(p-1)\int_Mw^{p/2}u_t\,\dd V.
			\label{ms:ent:Ivariation}
		\end{align}
		
		The divergence form of the $p$-Bochner identity is
		\begin{equation}\label{ms:ent:flux-Bochner}
		\begin{aligned}
			\Div\left(\frac12w^{p-2}A\nabla w
			-w^{p/2-1}(\Delta_pv)\nabla v\right)
			=w^{p-2}\bigl(|\nabla\nabla v|_A^2
			+\Ric(\nabla v,\nabla v)\bigr)-(\Delta_pv)^2.
		\end{aligned}
		\end{equation}
		Since the pressure flux is $C^1$, this identity
		holds in the sense of distributions across the critical set, with
		the continuous-extension convention above. Moreover,
		\begin{equation}\label{ms:ent:duality}
			(p-1)uw^{p/2-1}\langle\nabla v,\nabla w^{p/2}\rangle
			=\frac p2w^{p-2}\langle\nabla P,A\nabla w\rangle.
		\end{equation}
		This follows from $\nabla P=u\nabla v$,
		$A\nabla v=(p-1)\nabla v$, and the chain rule on $\{w>0\}$,
		and extends to the critical set under the same convention.
		
		Multiplying \eqref{ms:ent:flux-Bochner} by $P$ and
		integrating by parts therefore yields
		\begin{align*}
			-(p-1)\int_Mw^{p/2}u_t\,\dd V
			&=(p-1)\int_Muw^{p/2-1}
			\langle\nabla v,\nabla w^{p/2}\rangle\,\dd V\\
			&=\frac p2\int_Mw^{p-2}
			\langle\nabla P,A\nabla w\rangle\,\dd V\\
			&=p\int_Muw^{p/2}\Delta_pv\,\dd V-p\int_M P\left[
			w^{p-2}\bigl(|\nabla\nabla v|_A^2
			+\Ric(\nabla v,\nabla v)\bigr)
			-(\Delta_pv)^2\right]\,\dd V.
		\end{align*}
		Substitution into \eqref{ms:ent:Ivariation} proves
		\eqref{ms:ent:Iprime}.
	\end{proof}
	
	For the parameter $a>0$ in the main theorems, recall
	\begin{equation}\label{ms:ent:scalar-defect}
		D_a(s)
		=
		sP'(s)-\left(1-\frac1n\right)P(s)
		-\frac{a+1}{pa}\frac{P(s)^2}{U(s)}.
	\end{equation}
	By the normalization in Theorem~\ref{ms:thm:entropy},
	\begin{equation}\label{ms:ent:structural-condition}
		U>0,\qquad P>0
	\end{equation}
	on the solution range. We use the entropy functional
	\begin{equation}\label{ms:ent:NW-definition}
		\mathcal W_{F,p,a}
		=
		-\frac{\dd}{\dd t}\bigl(t^{a+1}\mathcal E\bigr)
		=
		t^{a+1}\mathcal I-(a+1)t^a\mathcal E.
	\end{equation}
	
	\begin{proof}[\bf Proof of Theorem~\ref{ms:thm:entropy}]
		The preceding variations give
		\begin{equation}\label{ms:ent:Wvariation}
			\mathcal W_{F,p,a}'
			=
			-t^{a+1}\mathcal E''
			+2(a+1)t^a\mathcal I
			-a(a+1)t^{a-1}\mathcal E,
		\end{equation}
		where
		\begin{equation}\label{ms:ent:Esecond}
			\begin{split}
				\mathcal E''
				=p\int_M\Bigg\{&
				P\left[
				w^{p-2}\bigl(|\nabla\nabla v|_A^2
				+\Ric(\nabla v,\nabla v)\bigr)
				-\frac{(\Delta_pv)^2}{n}
				\right]+
				\left[
				uP'(u)-\left(1-\frac1n\right)P
				\right](\Delta_pv)^2
				\Bigg\}\,\dd V.
			\end{split}
		\end{equation}
				We write the Hessian part of \eqref{ms:ent:Esecond} as a
		complete square. On the noncritical set, with $a_{ij}$ and
		$|\cdot|_A$ as in \eqref{ms:intro:WC-formula}, one has
		$A^{ij}a_{ij}=n$, $|a|_A^2=n$ and
		$w^{p/2-1}A^{ij}\nabla_i\nabla_jv=\Delta_pv$. Hence
		\begin{equation}\label{ms:ent:tracefree}
			w^{p-2}|\nabla\nabla v|_A^2
			-\frac{(\Delta_pv)^2}{n}
			=\left|w^{p/2-1}\nabla_i\nabla_jv-\frac{\Delta_pv}{n}a_{ij}\right|_A^2\ge0.
		\end{equation}
		Thus the Hessian part of \eqref{ms:ent:Esecond} is the $A$-norm of
		the trace-free part of $w^{p/2-1}\nabla\nabla v$; for $p=2$
		it is $|\nabla\nabla v-\frac{\Delta v}{n}g|^2$, the square in Ni's
		formula \cite{NiEntropy}. The identity extends to critical points
		under the continuous-extension convention above.
		
		For the scalar terms, completing the square gives
		\begin{align*}
			&\frac{a+1}{pa}\frac{P^2}{U}
			\left(\Delta_pv+\frac{aU}{tP}\right)^2=
			\frac{a+1}{pa}\frac{P^2}{U}(\Delta_pv)^2
			+\frac{2(a+1)}{pt}P\Delta_pv
			+\frac{a(a+1)}{pt^2}U.
		\end{align*}
				Combining this identity with
		\eqref{ms:ent:Eprime},
		\eqref{ms:ent:Wvariation},
		\eqref{ms:ent:Esecond}, and
		\eqref{ms:ent:tracefree}, we obtain
		\begin{equation}\label{ms:ent:W-dissipation}
			\begin{split}
				\mathcal W_{F,p,a}'
								=-pt^{a+1}\int_M\Bigg\{&
				P\left[
				\left|w^{p/2-1}\nabla_i\nabla_jv-\frac{\Delta_pv}{n}a_{ij}\right|_A^2
				+w^{p-2}\Ric(\nabla v,\nabla v)
				\right]\\
				&+D_a(u)(\Delta_pv)^2+\frac{a+1}{pa}\frac{P^2}{U}
				\left(\Delta_pv+\frac{aU}{tP}\right)^2
				\Bigg\}\,\dd V.
			\end{split}
		\end{equation}
		To pass to \eqref{ms:ent:Wsquare}, set $\kappa=aU/(tP)$. Since
		the trace-free part in \eqref{ms:ent:tracefree} is $A$-orthogonal to
		$a_{ij}$,
		\[
		\left|w^{p/2-1}\nabla_i\nabla_jv+\frac{\kappa}{n}a_{ij}\right|_A^2
		=\left|w^{p/2-1}\nabla_i\nabla_jv-\frac{\Delta_pv}{n}a_{ij}\right|_A^2
		+\frac1n\bigl(\Delta_pv+\kappa\bigr)^2,
		\]
		and, by \eqref{ms:ent:scalar-defect},
		\[
		\frac{a+1}{pa}\frac{P^2}{U}-\frac Pn=uP'(u)-P-D_a(u).
		\]
		Substituting these two identities into \eqref{ms:ent:W-dissipation}
		gives \eqref{ms:ent:Wsquare}.
		
		If $D_a\ge0$ on $\mathcal U$, every term in the integral in
		\eqref{ms:ent:W-dissipation} is nonnegative, which proves the
		monotonicity.
		If equality held at a positive time, the last square would vanish,
		so $\Delta_pv=-aU(u)/(tP(u))<0$ everywhere.
		This contradicts $\int_M\Delta_pv\,\dd V=0$.
		Thus the dissipation is strict on a closed manifold.
	\end{proof}
	
	\subsection{Entropy-power concavity}

The same second-variation identity yields concavity of
	the entropy power
	\begin{equation}\label{ms:ent:power-definition-body}
		\mathcal N_{F,p,a}(t)=\mathcal E(t)^{-1/a}.
	\end{equation}
	
	\begin{proof}[\bf Proof of Theorem~\ref{ms:thm:entropy-power}]
		Since $\mathcal E'=-\mathcal I$ and $\mathcal E>0$,
		differentiation gives
		\begin{equation}\label{ms:ent:power-first}
			\mathcal N_{F,p,a}'
			=
			\frac1a\mathcal E^{-1/a-1}\mathcal I
			\ge0,
		\end{equation}
		and
		\begin{equation}\label{ms:ent:power-chain}
			\mathcal N_{F,p,a}''
			=
			-\frac1a\mathcal E^{-1/a-1}
			\left[
			\mathcal E''
			-\frac{a+1}{a}\frac{\mathcal I^2}{\mathcal E}
			\right].
		\end{equation}
		On the other hand,
		$\int_M P\Delta_pv\,\dd V=-\mathcal I$ gives
		\begin{equation}\label{ms:ent:power-variance}
			\begin{split}
				&\int_M U
				\left(
				\frac{P\Delta_pv}{U}
				+\frac{\mathcal I}{\mathcal E}
				\right)^2\,\dd V=
				\int_M\frac{P^2}{U}(\Delta_pv)^2\,\dd V
				-\frac{\mathcal I^2}{\mathcal E}.
			\end{split}
		\end{equation}
		Substituting \eqref{ms:ent:Esecond} and
		\eqref{ms:ent:power-variance} into
		\eqref{ms:ent:power-chain} yields
		\begin{equation}\label{ms:ent:power-dissipation}
			\begin{split}
				\mathcal N_{F,p,a}''
				=-\frac pa\mathcal E^{-1/a-1}
				\int_M\Bigg\{&
				P\left[
				w^{p-2}\bigl(|\nabla\nabla v|_A^2
				+\Ric(\nabla v,\nabla v)\bigr)
				-\frac{(\Delta_pv)^2}{n}
				\right]\\
				&+D_a(u)(\Delta_pv)^2+\frac{a+1}{pa}U
				\left(
				\frac{P\Delta_pv}{U}
				+\frac{\mathcal I}{\mathcal E}
				\right)^2
				\Bigg\}\,\dd V
				\le0.
			\end{split}
		\end{equation}
		If $\mathcal N_{F,p,a}''=0$ at some time, the last square vanishes.
		Integrating
		$\Delta_pv=-(\mathcal I/\mathcal E)U(u)/P(u)$ over $M$ gives
		\[
		0=-\frac{\mathcal I}{\mathcal E}
		\int_M\frac{U(u)}{P(u)}\,\dd V.
		\]
		Since $U,P,\mathcal E>0$, this forces $\mathcal I=0$.
		Hence $v$, and therefore $u$, is constant on each connected
		component. Conversely, such a density makes every term in
		\eqref{ms:ent:power-dissipation} zero.
		If the value on a component is $k$ at time $t_0$, then
		\[
		\frac{\dd}{\dd t}\int (u-k)^2\,\dd V
		=-2\int F'(u)^{p-1}|\nabla u|^p\,\dd V\le0
		\]
		on that component. The zero initial value proves stationarity
		for subsequent times.
		
		Finally, comparing \eqref{ms:ent:Wvariation} with
		\eqref{ms:ent:power-chain} gives
		\begin{equation}\label{ms:ent:W-power-identity}
			\begin{split}
				\mathcal W_{F,p,a}'
				={}&
				at^{a+1}\mathcal E^{1+1/a}
				\mathcal N_{F,p,a}''-\frac{a+1}{a}t^{a+1}\mathcal E
				\left(
				\frac{\mathcal I}{\mathcal E}
				-\frac at
				\right)^2.
			\end{split}
		\end{equation}
		This identity relates entropy-power concavity to
		$W$-entropy dissipation.
	\end{proof}
	
	We next consider the pointwise equality system, which
	determines both the diffusion law on the attained
	density interval and the local geometry.
	
	\begin{proof}[Proof of Theorem~\ref{ms:thm:rigidity}]
		Work at the fixed time and on the connected open set
		$\Omega$ in the theorem.
		Since $|\nabla v|>0$, equality in
		\eqref{ms:ent:tracefree}, together with
		\eqref{ms:ent:critical-square}, gives
		\begin{equation}\label{ms:rig:local-flux-trace}
			\nabla\bigl(|\nabla v|^{p-2}\nabla v\bigr)=\frac{\Delta_pv}{n}\operatorname{id}.
		\end{equation}
		Commuting covariant derivatives and contracting,
		we obtain
		\begin{equation}\label{ms:rig:contracted-identity}
			(n-1)\nabla\left(\frac{\Delta_pv}{n}\right)
			=-|\nabla v|^{p-2}\Ric(\nabla v,\cdot).
		\end{equation}
		
		Since $\Ric\ge0$ and $\Ric(\nabla v,\nabla v)=0$, one has
		$\Ric(\nabla v,\cdot)=0$.
		Thus $\Delta_pv$ is constant on $\Omega$.
		
		The last equality in
		\eqref{ms:intro:local-equality-system} gives
		\[
		\Delta_pv=-\lambda\frac{U}{P}<0.
		\]
		Hence
		\begin{equation}\label{ms:rig:local-ratio}
			\frac{P(u)}{U(u)}
			=-\frac{\lambda}{\Delta_pv}
			=:\delta>0
		\end{equation}
		is constant on $\Omega$.
		Because $U''>0$ and $|\nabla v|>0$, the density has
		no critical points on $\Omega$, and $u(\Omega)$
		is an open interval.
		On this interval, $P=\delta U$ implies
		\[
		sU'(s)=(1+\delta)U(s).
		\]
		Therefore
		\begin{equation}\label{ms:rig:local-energy-law}
			U(s)=C_Us^{1+\delta},
			\qquad C_U>0.
		\end{equation}
		In particular,
		\[
		sP'(s)=(1+\delta)P(s),\qquad
		\frac{P(s)^2}{U(s)}=\delta P(s).
		\]
		Substituting these identities into $D_a=0$ gives
		\[
		\delta+\frac1n
		=\frac{a+1}{pa}\delta.
		\]
		Since $\delta>0$, it follows that
		\begin{equation}\label{ms:rig:local-parameters}
			0<a<\frac1{p-1},\qquad
			\delta=\frac{pa}{n[1-(p-1)a]}.
		\end{equation}
		Combining \eqref{ms:rig:local-flux-trace} and
		\eqref{ms:rig:local-ratio}, we obtain
		\begin{equation}\label{ms:rig:local-flux}
			\nabla\bigl(|\nabla v|^{p-2}\nabla v\bigr)
			=-\frac{\lambda}{n\delta}\operatorname{id}.
		\end{equation}
		
		By \eqref{ms:ent:Upp} and
		\eqref{ms:rig:local-energy-law},
		\[
		F'(s)
		=s^{1/(p-1)}U''(s)
		=C_U\delta(1+\delta)
		s^{\delta-1+1/(p-1)}.
		\]
		Integrating on $u(\Omega)$ yields
		\begin{equation}\label{ms:rig:local-power-law}
			F(s)=C_Fs^{\delta+1/(p-1)}+C_0,
			\qquad
			C_F=
			\frac{C_U\delta(1+\delta)}
			{\delta+1/(p-1)}>0.
		\end{equation}
		
		To identify the local metric, set
		\[
		r=\frac{n\delta}{\lambda}|\nabla v|^{p-1}.
		\]
		Equation \eqref{ms:rig:local-flux} gives
		\[
		\nabla r=-\frac{\nabla v}{|\nabla v|},\qquad
		|\nabla r|=1,\qquad
		\nabla\nabla r
		=\frac1r\bigl(g-\dd r\otimes\dd r\bigr).
		\]
		Consequently, in local coordinates adapted to
		the level sets of $r$,
		\begin{equation}\label{ms:rig:local-cone}
			g=\dd r^2+r^2g_\Sigma,\qquad
			|\nabla v|^{p-2}\nabla v=-\frac{\lambda}{n\delta}r\,\partial_r,
		\end{equation}
		where $g_\Sigma$ is independent of $r$.
		This proves the geometric conclusion.
		
		The corresponding pressure profile follows from
		the inverse flux relation:
		\[
		\nabla v
		=
		-\left(\frac{\lambda}{n\delta}r\right)^{1/(p-1)}
		\partial_r.
		\]
		Hence, locally,
		\[
		v
		=
		v_0-\frac{p-1}{p}
		\left(\frac{\lambda}{n\delta}\right)^{1/(p-1)}
		r^{p/(p-1)}.
		\]
	\end{proof}

	\section{Parabolic frequency monotonicity}
	\label{ms:freq:section}
	
	Throughout this section, $(M,g)$ is closed, with fixed metric
	and volume element $\dd V$. The frequency formula requires
	no curvature assumption.
	
	Fix a reference density $\ell$. Assume that $F$ is continuous
	and strictly increasing on an interval containing $\ell$ and
	the solution range, and that $F\in C^2$ with $F'>0$ near the
	attained values. Define
	\begin{equation}\label{ms:freq:primitive}
		\Psi_\ell(s)=\int_\ell^s[F(r)-F(\ell)]\,\dd r,
	\end{equation}
	and
	\begin{equation}\label{ms:freq:adapted-def}
		I(t)=2\int_M\Psi_\ell(u)\,\dd V,\qquad
		D(t)=-\int_M|\nabla F(u)|^p\,\dd V.
	\end{equation}
	Strict monotonicity of $F$ gives $I\ge0$, with equality
	precisely when $u=\ell$ almost everywhere, while $D\le0$.
	
	We first carry out the calculations for positive solutions of
	\eqref{ms:intro:PF} satisfying
	\begin{equation}\label{ms:freq:classical-class}
		u\in C^1_{\mathrm{loc}}\bigl((0,T);C^1(M)\bigr),
	\end{equation}
	where the equation holds in the sense of distributions on $M$
	for each $t$.
	
	\begin{lemma}[Energy identities]
		\label{ms:freq:energy-lemma}
		The energies in \eqref{ms:freq:adapted-def} satisfy
		\begin{equation}\label{ms:freq:energy-identities}
			I'=2D,\qquad
			D'=p\int_MF'(u)u_t^2\,\dd V.
		\end{equation}
		Moreover,
		\begin{equation}\label{ms:freq:mass}
			\frac{\dd}{\dd t}\int_Mu\,\dd V=0,
		\end{equation}
		and
		\begin{equation}\label{ms:freq:pairing}
			D=\int_M[F(u)-F(\ell)]u_t\,\dd V.
		\end{equation}
	\end{lemma}
	
	\begin{proof}
		Integrating the equation over $M$ proves
		\eqref{ms:freq:mass}. Integration by parts gives
		\begin{align*}
			I'
			&=2\int_M[F(u)-F(\ell)]u_t\,\dd V\\
			&=2\int_M[F(u)-F(\ell)]\Delta_pF(u)\,\dd V\\
			&=-2\int_M|\nabla F(u)|^p\,\dd V
			=2D,
		\end{align*}
		which also proves \eqref{ms:freq:pairing}.
		
		Since $\xi\mapsto|\xi|^p$ is $C^1$ for every $p>1$,
		differentiation of the gradient energy yields
		\begin{align*}
			D'
			&=-p\int_M|\nabla F(u)|^{p-2}
			\left\langle\nabla F(u),
			\nabla\bigl(F'(u)u_t\bigr)\right\rangle\,\dd V\\
			&=p\int_MF'(u)u_t\Delta_pF(u)\,\dd V\\
			&=p\int_MF'(u)u_t^2\,\dd V.
		\end{align*}
		These integrations use the continuous flux
		$|\nabla F(u)|^{p-2}\nabla F(u)$, including at its zeros.
	\end{proof}
	
	The same identities hold almost everywhere for energy solutions
	satisfying, on each compact positive-time interval $J$,
	\begin{equation}\label{ms:freq:weak-class}
		\begin{gathered}
			0<u_-\le u\le u_+,\qquad
			u\in W^{1,2}(J;L^2(M)),\\
			F(u)\in L^p(J;W^{1,p}(M)),\qquad
			u_t=\Delta_pF(u)\quad\text{in distributions},
		\end{gathered}
	\end{equation}
	where $F\in C^2$ and $F'>0$ near $[u_-,u_+]$.
	Indeed, $(F(u))_t=F'(u)u_t$, and the convex chain rule for
	the $p$-Dirichlet energy gives
	\[
	\frac{\dd}{\dd t}
	\left(\frac1p\int_M|\nabla F(u)|^p\,\dd V\right)
	=-\int_MF'(u)u_t^2\,\dd V.
	\]
	The chain rule for $\int_M\Psi_\ell(u)\,\dd V$ gives the
	first identity. In this class, $I$ and $D$ are taken with their
	locally absolutely continuous representatives.
	
	Suppose that a fixed constant $\kappa>0$ satisfies
	\begin{equation}\label{ms:freq:scalar-condition}
		[F(s)-F(\ell)]^2
		\le\kappa F'(s)\Psi_\ell(s)
	\end{equation}
	on the attained density range. On each interval where $I>0$,
	define
	\begin{equation}\label{ms:freq:quotient-definition}
		N_\gamma(t)=\frac{D(t)}{I(t)^\gamma},
		\qquad 0<\gamma\le\frac p\kappa.
	\end{equation}
	
	\begin{proof}[\bf Proof of Theorem \ref{ms:thm:frequency} ]
		By \eqref{ms:freq:energy-identities},
		\begin{equation}\label{ms:freq:quotient}
			N_\gamma'
			=I^{-\gamma}
			\left[
			p\int_MF'(u)u_t^2\,\dd V
			-\frac{2\gamma D^2}{I}
			\right].
		\end{equation}
		Weighted Cauchy--Schwarz and
		\eqref{ms:freq:scalar-condition} imply
		\begin{align}
			D^2
			&=
			\left(\int_M[F(u)-F(\ell)]u_t\,\dd V\right)^2
			\notag\\
			&\le
			\left(\int_M\frac{[F(u)-F(\ell)]^2}{F'(u)}
			\,\dd V\right)
			\left(\int_MF'(u)u_t^2\,\dd V\right)
			\notag\\
			&\le\frac{\kappa I}{2}
			\int_MF'(u)u_t^2\,\dd V.
			\label{ms:freq:weighted-CS}
		\end{align}
		Consequently,
		\begin{equation}\label{ms:freq:derivative}
			N_\gamma'
			\ge(p-\gamma\kappa)I^{-\gamma}
			\int_MF'(u)u_t^2\,\dd V
			\ge0.
		\end{equation}
		
		To retain the equality information, expand the square
		with coefficient $2D/(\kappa I)$. This gives the exact
		decomposition
		\begin{equation}\label{ms:freq:three-deficits}
			\begin{aligned}
				N_\gamma'={}&
				\frac p{I^\gamma}
				\int_MF'(u)
				\left[
				u_t-\frac{2D}{\kappa I}
				\frac{F(u)-F(\ell)}{F'(u)}
				\right]^2\,\dd V+\frac{4pD^2}{\kappa^2I^{\gamma+2}}
				\int_M
				\left[
				\kappa\Psi_\ell(u)
				-\frac{[F(u)-F(\ell)]^2}{F'(u)}
				\right]\,\dd V\\
				&+2\left(\frac p\kappa-\gamma\right)
				\frac{D^2}{I^{\gamma+1}}.
			\end{aligned}
		\end{equation}
		Indeed, the first integral expands as
		\[
		\int_MF'(u)u_t^2\,\dd V
		-\frac{4D^2}{\kappa I}
		+\frac{4D^2}{\kappa^2I^2}
		\int_M\frac{[F(u)-F(\ell)]^2}{F'(u)}\,\dd V,
		\]
		and substitution into \eqref{ms:freq:three-deficits}
		recovers \eqref{ms:freq:quotient}.
		
		Every term in \eqref{ms:freq:three-deficits} is
		nonnegative. For energy solutions, the calculation holds
		almost everywhere, and local absolute continuity proves
		the same monotonicity.
	\end{proof}
	
	The scalar condition \eqref{ms:freq:structure} involves only the
	density range. If $F\in C^2$ and $F'>0$ on a compact interval containing
	$\ell$, then
	\begin{equation}\label{ms:freq:ratio-limit}
		\lim_{s\to\ell}
		\frac{[F(s)-F(\ell)]^2}
		{F'(s)\Psi_\ell(s)}=2.
	\end{equation}
	Thus the ratio extends continuously at $\ell$, and its maximum
	on that interval provides a finite admissible $\kappa$.
	
	\subsection{Time convexity and nonvanishing}

The frequency also determines the time convexity of the energy.
	On an interval where $I>0$, set
	\begin{equation}\label{ms:freq:time-transform}
		\mathcal J_\gamma(t)=
		\begin{cases}
			\displaystyle
			\frac{I(t)^{1-\gamma}-1}{1-\gamma},
			&\gamma\ne1,\\[2mm]
			\log I(t),&\gamma=1.
		\end{cases}
	\end{equation}
	Then
	\begin{equation}\label{ms:freq:time-convexity}
		\mathcal J_\gamma'=2N_\gamma,\qquad
		\mathcal J_\gamma''=2N_\gamma'\ge0,
	\end{equation}
	with the second identity interpreted distributionally for
	energy solutions. In particular, for $t_1<t_2<t_3$,
	\begin{equation}\label{ms:freq:three-time}
		\mathcal J_\gamma(t_2)
		\le
		\frac{t_3-t_2}{t_3-t_1}\mathcal J_\gamma(t_1)
		+\frac{t_2-t_1}{t_3-t_1}\mathcal J_\gamma(t_3).
	\end{equation}
	
	If an admissible $\gamma\ge1$ is available, the supporting-line
	inequality
	\[
	\mathcal J_\gamma(t)
	\ge\mathcal J_\gamma(s)+2N_\gamma(s)(t-s)
	\]
	gives, for $t\ge s$,
	\begin{equation}\label{ms:freq:nonvanishing}
		I(t)\ge
		\begin{cases}
			I(s)\exp\!\bigl(2N_1(s)(t-s)\bigr),
			&\gamma=1,\\[2mm]
			\left[
			I(s)^{1-\gamma}
			-2(\gamma-1)N_\gamma(s)(t-s)
			\right]^{-1/(\gamma-1)},
			&\gamma>1.
		\end{cases}
	\end{equation}
	The expressions on the right are positive whenever $I(s)>0$,
	since $N_\gamma\le0$. Applying the estimate up to a hypothetical
	first zero of $I$ excludes arrival at the reference constant in
	finite time under this exponent condition. If $I(s)=0$,
	the identity $I'=2D\le0$ implies $I(t)=0$ for all later times.

\section{Examples}
\label{ms:examples:section}

\subsection{The case \texorpdfstring{$p=2$}{p=2}}
\label{ms:examples:p2}

For $p=2$, equation \eqref{ms:intro:PF} is the filtration equation
$u_t=\Delta F(u)$ studied in \cite{Yau94,Ma,Xu,HaoWang}. In this case
$A=\operatorname{id}$, $\Delta_p=\Delta$, and the pressure
\eqref{ms:intro:pressure} is $v=G(u)$, where $G'(s)=F'(s)/s$ is the
function in \eqref{ms:intro:Yau}. Moreover, $\psi'=F'(s)$ and $b(s)=sF''(s)/F'(s)$. We first check the
hypotheses of the main theorems.

\emph{Li--Yau estimate.} We do this for general $p\ge2$, for later use.
For $\alpha=p$, one has $(p-1)sF''+(p-2)F'=(p-1)bF'$, and
\eqref{ms:eq:LY-polynomial} becomes
\begin{equation}\label{ms:examples:LY-reduced}
	-\left(b+\frac1n\right)
	\left[\frac{(p-1)^2}{p}\psi'y^2+\frac{2(p-1)}{p}c(t)y
	+\frac{c(t)^2}{p\psi'}\right]-c'(t)\le0.
\end{equation}
Suppose that $\mu:=\min_{\mathcal U}b>-1/n$, and let $c(t)=K/t$ with
$K=p\max_{\mathcal U}\psi'/(\mu+1/n)$. The terms in $y^2$ and $y$ are
nonpositive, and since $c'=-c^2/K$, the remaining terms are at most
$c^2[1/K-(\mu+1/n)/(p\max_{\mathcal U}\psi')]=0$. Applying
Theorem~\ref{ms:thm:LY} to $u(\cdot,\cdot+\tau)$ with $c(s)=K/s$ from a
small initial time $\varepsilon$, where $K/\varepsilon$ dominates
$|\nabla v|^p-pv_t$, and letting $\varepsilon,\tau\to0$, we obtain
\begin{equation}\label{ms:examples:LY-explicit}
	|\nabla v|^p-pv_t\le\frac{p\max_{\mathcal U}\psi'}{(\mu+1/n)\,t}
	\qquad\text{on }M\times(0,T).
\end{equation}

\emph{Hamilton estimate.} For $p=2$, condition
\eqref{ms:eq:Hamilton-condition} reads
\begin{equation}\label{ms:examples:Hamilton-p2}
	q+qb(s)+q\left(1-\frac q2\right)\frac{F'(s)}{\alpha-G(s)}
	-\frac{(n-1)[\alpha-G(s)]}{2F'(s)}b(s)^2\ge\beta,
	\qquad s\in\mathcal U.
\end{equation}

\emph{Entropy.} Suppose that $\int_0^1r^{-1}F'(r)\,\dd r<\infty$. Then
$G(0)$ and $F(0)=\lim_{s\to0}F(s)$ are finite. With $G(0)=0$,
\eqref{ms:intro:energy-notation} and \eqref{ms:intro:scalar-defect}
become
\begin{equation}\label{ms:examples:entropy-p2}
	\begin{gathered}
		U(s)=\int_0^sG(r)\,\dd r,\qquad
		P(s)=F(s)-F(0),\qquad
		sP'(s)=sF'(s),\\
		D_a(s)=sF'(s)-\left(1-\frac1n\right)[F(s)-F(0)]
		-\frac{a+1}{2a}\frac{[F(s)-F(0)]^2}{U(s)}.
	\end{gathered}
\end{equation}

The main theorems therefore give the following.

\begin{corollary}[The case $p=2$]\label{ms:cor:p2}
	Let $F\in C^3$ with $F'>0$ on a neighborhood of $\mathcal U$, let
	$G'(s)=F'(s)/s$, and let $u>0$ be a smooth solution of
	$u_t=\Delta F(u)$ with values in $\mathcal U$.
	\begin{enumerate}[label=\rm(\roman*)]
		\item \emph{Li--Yau estimate (Theorem~\ref{ms:thm:LY}).}
		If $M$ is closed with $\Ric\ge0$ and
		$\mu:=\min_{s\in\mathcal U}sF''(s)/F'(s)>-1/n$, then
		\[
		\frac{F'(u)^2|\nabla u|^2}{u^2}-\frac{2F'(u)u_t}{u}
		=|\nabla G(u)|^2-2\partial_tG(u)
		\le\frac{2\max_{\mathcal U}F'}{(\mu+1/n)\,t}
		\qquad\text{on }M\times(0,T).
		\]
		\item \emph{Aronson--B\'enilan estimate (Theorem~\ref{ms:thm:AB}).}
		If $M$ is closed with $\Ric\ge0$, $sF''(s)/F'(s)$ is nonincreasing
		on $\mathcal U$, and
		$\varrho:=u_+F''(u_+)/F'(u_+)+2/n>0$, then
		\[
		\Delta G(u)\ge-\frac1{\varrho t},
		\qquad\text{equivalently}\qquad
		\frac{u_t}{u}\ge F'(u)\frac{|\nabla u|^2}{u^2}-\frac1{\varrho t}
		\qquad\text{on }M\times(0,T).
		\]
		\item \emph{Hamilton-type estimate (Theorem~\ref{ms:thm:Hamilton}).}
		Let $M$ be complete with $\Ric\ge-\kappa g$ on $B(x_0,R)$
		{\rm(}which suffices for $p=2$; see the end of the proof of
		Theorem~\ref{ms:thm:Hamilton}{\rm)}, let $u$
		solve the equation in $Q_{R,T}$, and suppose that
		\eqref{ms:examples:Hamilton-p2} holds for some
		$\alpha>\max_{\mathcal U}G$ and $q,\beta>0$. Then
		\[
		\sup_{Q_{R/2,T/2}}
		\frac{F'(u)|\nabla u|}{u\,[\alpha-G(u)]^{q/2}}
		\le C\left(\frac1R+\frac1{\sqrt T}+\sqrt\kappa\right),
		\]
		with $C$ as in Theorem~\ref{ms:thm:Hamilton}.
		\item \emph{$\mathcal W$-entropy formula
			(Theorem~\ref{ms:thm:entropy}).}
		Let $M$ be closed with $\Ric\ge0$, suppose that
		$\int_0^1r^{-1}F'(r)\,\dd r<\infty$, and let $U$ and $D_a$ be given by
		\eqref{ms:examples:entropy-p2}. If $D_a\ge0$ on $\mathcal U$ for
		some $a>0$, then
		\[
		\mathcal W_{F,2,a}(t)
		=t^{a+1}\int_M\left[\frac{F'(u)^2|\nabla u|^2}{u}
		-\frac{a+1}{t}U(u)\right]\dd V
		\]
		is strictly decreasing on $(0,T)$, and
		\[
		\begin{aligned}
			\frac{\dd}{\dd t}\mathcal W_{F,2,a}
			={}&-2t^{a+1}\int_M[F(u)-F(0)]
			\left(\left|\nabla\nabla G(u)
			+\frac{aU(u)}{nt[F(u)-F(0)]}g\right|^2
			+\Ric\bigl(\nabla G(u),\nabla G(u)\bigr)\right)\dd V\\
			&-2t^{a+1}\int_M\bigl[uF'(u)-F(u)+F(0)-D_a(u)\bigr]
			\left(\Delta G(u)+\frac{aU(u)}{t[F(u)-F(0)]}\right)^2\dd V\\
			&-2t^{a+1}\int_MD_a(u)\bigl(\Delta G(u)\bigr)^2\,\dd V.
		\end{aligned}
		\]
		\item \emph{Entropy-power concavity
			(Theorem~\ref{ms:thm:entropy-power}).}
		Under the hypotheses of {\rm(iv)},
		$\mathcal N_{F,2,a}(t)=\bigl(\int_MU(u)\,\dd V\bigr)^{-1/a}$ is
		nondecreasing and concave on $(0,T)$, and
		$\mathcal N_{F,2,a}''(t_0)=0$ if and only if $u(\cdot,t_0)$ is
		constant on each connected component of $M$.
	\end{enumerate}
\end{corollary}

\begin{remark}\label{ms:rem:p2}
	For $F(s)=s^m$, part~(ii) is \eqref{ms:intro:LNVV-AB}, and parts
	(iv)--(v) with $m>1$ and $a=n(m-1)/(n(m-1)+2)$ give the formula of
	\cite{LNVV} and the R\'enyi entropy power \eqref{ms:intro:Renyi},
	since $U(s)=s^m/(m-1)$ and $D_a\equiv0$. For a general $F$,
	Hao--Wang \cite{HaoWang} derived local differential Harnack estimates
	by Nash--Moser iteration on complete manifolds with Ricci curvature
	bounded below, together with the frequency monotonicity that
	Theorem~\ref{ms:thm:frequency} recovers for $p=2$. Parts (i)--(iii)
	complement these with maximum-principle estimates: (i) is an explicit
	instance of Yau's estimate \eqref{ms:intro:Yau}, (ii) extends
	\eqref{ms:intro:LNVV-AB} to every $F$ with $sF''/F'$ nonincreasing and \(\inf_{\mathcal U}sF''(s)/F'(s)>-2/n\), and
	(iii) is Xu's estimate \eqref{ms:intro:Xu} with a free exponent $q$.
	Parts~(iv) and (v) have no counterpart in \cite{HaoWang}; to our
	knowledge, they are the first $\mathcal W$-entropy formula and
	entropy-power concavity for the filtration equation with a general
	nonlinearity $F$, such formulas having been known only for
	$F(s)=s^m$ \cite{LNVV,ST,LLRenyi}.
\end{remark}

\subsection{Positive sums of powers}
\label{ms:examples:mixed-powers}

Let
\[
F(s)=\sum_{j=1}^N c_js^{m_j},
\qquad c_j>0,\qquad 0<m_1<\cdots<m_N,\qquad N\ge2.
\]
For $p=2$, $u_t=\Delta F(u)$ is a generalized porous-medium equation
whose diffusivity behaves like different powers at small and large
densities \cite{VazquezPME,FotacheMuratori}; when all but one
coefficient vanish, it reduces to the porous-medium or fast-diffusion
equation.

Put $\delta_j=m_j-1/(p-1)$. Then
\begin{equation}\label{ms:examples:sum-pressure}
	\psi'=\sum_jc_jm_js^{\delta_j},
	\qquad
	b=
	\frac{\sum_j\delta_jc_jm_js^{\delta_j}}
	{\sum_jc_jm_js^{\delta_j}},
	\qquad
	sb'(s)=
	\frac{\sum_j\delta_j^2c_jm_js^{\delta_j}}
	{\sum_jc_jm_js^{\delta_j}}-b^2>0.
\end{equation}
Thus $b$ is a weighted mean of the $\delta_j$, and $sb'$ is the
corresponding variance. Hence $b$ is strictly increasing,
$\delta_1<b<\delta_N$, and $\min_{\mathcal U}b=b(u_-)$.
We now check the hypotheses of the main theorems.

\emph{Li--Yau estimate.} By Section~\ref{ms:examples:p2},
\eqref{ms:examples:LY-explicit} holds when $\mu=b(u_-)>-1/n$, in
particular when $\delta_1\ge-1/n$.

\emph{Aronson--B\'enilan estimate.} Since $b'>0$, the monotonicity
hypothesis of Theorem~\ref{ms:thm:AB} fails.

\emph{Hamilton estimate.} Let
\begin{equation}\label{ms:examples:Hamilton-q}
	0<q<\frac{p(p-1)}{(p-1)^2+(n-1)(p-2)^2/4}.
\end{equation}
Then the coefficient of $\psi'/(\alpha-v)$ in
\eqref{ms:eq:Hamilton-condition} is positive. Fix $u_0>0$, let
$J=[u_0/2,2u_0]$, and let $\mathcal U\subset J$ and
$\alpha=\max_{\mathcal U}v+\varepsilon$. Then
$\psi'/(\alpha-v)\ge\min_J\psi'/(\operatorname{osc}_{\mathcal U}v+\varepsilon)$
on $\mathcal U$, and since $b$ and $\psi'$ are bounded on $J$, the
left-hand side of \eqref{ms:eq:Hamilton-condition} is at least $q(p-1)$
once this lower bound is large. Hence there are $\eta,\varepsilon_0>0$, depending only on
$n,p,q,u_0$ and $F|_J$, such that \eqref{ms:eq:Hamilton-condition}
holds with $\beta=q(p-1)$ whenever
$\mathcal U\subset[u_0-\eta,u_0+\eta]$ and $0<\varepsilon\le\varepsilon_0$.
This argument uses only that $F\in C^3$ and $F'>0$ on $J$.

\emph{Entropy.} Assume that $\delta_j>0$ for every $j$. Then the
pressure integral converges at zero, and with $v(0)=0$,
\begin{equation}\label{ms:examples:sum-energies}
	\begin{gathered}
		v=\sum_j\frac{c_jm_j}{\delta_j}s^{\delta_j},
		\qquad
		U=\sum_jU_j,\qquad
		U_j(s)=\frac{c_jm_j}{\delta_j(\delta_j+1)}s^{\delta_j+1},\\
		P=\sum_j\delta_jU_j,
		\qquad
		sP'-P=\sum_j\delta_j^2U_j,
	\end{gathered}
\end{equation}
so that
\begin{equation}\label{ms:examples:sum-defect}
	D_a=\sum_j\delta_j^2U_j+\frac1n\sum_j\delta_jU_j
	-\frac{a+1}{pa}
	\frac{\bigl(\sum_j\delta_jU_j\bigr)^2}{\sum_jU_j}.
\end{equation}
Let $a_*=n\delta_N/(n(p-1)\delta_N+p)$. Then
$\frac{a_*+1}{pa_*}=1+\frac1{n\delta_N}$, and Lagrange's identity gives
\begin{equation}\label{ms:examples:mixed-defect}
	UD_{a_*}
	=\sum_{i<j}(\delta_i-\delta_j)^2U_iU_j
	+\frac Pn\left(U-\frac P{\delta_N}\right)>0,
\end{equation}
since $P\le\delta_NU$. As $(a+1)/(pa)$ is decreasing in $a$,
$D_a\ge D_{a_*}>0$ on $(0,\infty)$ for every $a\ge a_*$. The variance
term in \eqref{ms:examples:mixed-defect} vanishes only for a single
power; by real analyticity, it also excludes the equality case of
Theorem~\ref{ms:thm:rigidity}.

\emph{Frequency.} Take $\ell=0$. Then
$\Psi_0(s)=\sum_j\frac{c_j}{m_j+1}s^{m_j+1}$, and the Cauchy--Schwarz
inequality gives
\[
F(s)^2
\le F'(s)\sum_j\frac{c_j}{m_j}s^{m_j+1}
\le\left(1+\frac1{m_1}\right)F'(s)\Psi_0(s).
\]
Hence $\kappa=1+1/m_1$ is admissible in \eqref{ms:freq:structure}.

\begin{corollary}[Positive sums of powers]\label{ms:cor:sums}
	Let $F(s)=\sum_{j=1}^Nc_js^{m_j}$ be as above, let
	$\delta_j=m_j-1/(p-1)$, and let $u>0$ be a smooth solution of
	\eqref{ms:intro:PF} with values in $\mathcal U$.
	\begin{enumerate}[label=\rm(\roman*)]
		\item \emph{Li--Yau estimate (Theorem~\ref{ms:thm:LY}).}
		If $p\ge2$, $M$ is closed with $\Ric\ge0$, and
		\[
		\mu:=\frac{\sum_j\delta_jc_jm_ju_-^{\delta_j}}
		{\sum_jc_jm_ju_-^{\delta_j}}>-\frac1n,
		\]
		then, on $M\times(0,T)$,
		\[
		\frac{\bigl|\nabla\sum_jc_ju^{m_j}\bigr|^p}{u^{p/(p-1)}}
		-p\,\frac{\sum_jc_jm_ju^{m_j-1}}{u^{1/(p-1)}}\,u_t
		\le\frac{p\max_{s\in\mathcal U}\sum_jc_jm_js^{\delta_j}}
		{(\mu+1/n)\,t}.
		\]
		\item \emph{Hamilton-type estimate (Theorem~\ref{ms:thm:Hamilton}).}
		Let $p>1$, let $q$ satisfy \eqref{ms:examples:Hamilton-q}, let
		$u_0>0$, and let $\eta,\varepsilon_0$ be as above. Suppose that
		$\mathcal U\subset[u_0-\eta,u_0+\eta]$,
		$\alpha=\max_{\mathcal U}v+\varepsilon$ with
		$0<\varepsilon\le\varepsilon_0$, $M$ is complete with
		$\operatorname{Sec}\ge-\kappa$ on $B(x_0,R)$, and $u$ solves
		\eqref{ms:intro:PF} in $Q_{R,T}$. Then
		\[
		\sup_{Q_{R/2,T/2}}
		\frac{\bigl|\nabla\sum_jc_ju^{m_j}\bigr|}
		{u^{1/(p-1)}(\alpha-v)^{q/p}}
		\le C\left(\frac1R+\frac1{T^{1/p}}+\sqrt\kappa\right),
		\]
		where $v$ is the pressure \eqref{ms:intro:pressure} and $C$ depends
		only on $n,p,q,\varepsilon,u_0$ and the $c_j,m_j$.
		\item \emph{$\mathcal W$-entropy formula
			(Theorem~\ref{ms:thm:entropy}).}
		Suppose that $\delta_j>0$ for every $j$, that $a\ge a_*$, that $M$ is
		closed with $\Ric\ge0$, and that $u$ satisfies the regularity
		hypothesis of Theorem~\ref{ms:thm:entropy}. Let $U_j$ and $D_a$ be
		given by \eqref{ms:examples:sum-energies} and
		\eqref{ms:examples:sum-defect}. Then
		\[
		\mathcal W_{F,p,a}(t)
		=t^{a+1}\int_M\left[
		\frac{\bigl|\nabla\sum_jc_ju^{m_j}\bigr|^p}{u^{1/(p-1)}}
		-\frac{a+1}{t}\sum_{j=1}^N
		\frac{c_jm_j}{\delta_j(\delta_j+1)}u^{\delta_j+1}\right]\dd V
		\]
		is strictly decreasing on $(0,T)$, and
		\[
		\begin{aligned}
			\frac{\dd}{\dd t}\mathcal W_{F,p,a}
			={}&-pt^{a+1}\int_M
			\Bigl(\sum_j\delta_jU_j(u)\Bigr)
			\left(\left||\nabla v|^{p-2}\nabla_i\nabla_jv
			+\frac{a\sum_jU_j(u)}{nt\sum_j\delta_jU_j(u)}a_{ij}\right|_A^2
			+|\nabla v|^{2p-4}\Ric(\nabla v,\nabla v)\right)\dd V\\
			&-pt^{a+1}\int_M
			\Bigl[\sum_j\delta_j^2U_j(u)-D_a(u)\Bigr]
			\left(\Delta_pv+\frac{a\sum_jU_j(u)}{t\sum_j\delta_jU_j(u)}\right)^2\dd V\\
			&-pt^{a+1}\int_MD_a(u)(\Delta_pv)^2\,\dd V.
		\end{aligned}
		\]
		\item \emph{Entropy-power concavity
			(Theorem~\ref{ms:thm:entropy-power}).}
		Under the hypotheses of {\rm(iii)}, the entropy power
		$\mathcal N_{F,p,a}=\bigl(\int_M\sum_jU_j(u)\,\dd V\bigr)^{-1/a}$
		is nondecreasing and concave on $(0,T)$, and
		$\mathcal N_{F,p,a}''(t_0)=0$ if and only if $u(\cdot,t_0)$ is
		constant on each connected component of $M$.
		\item \emph{Frequency monotonicity (Theorem~\ref{ms:thm:frequency}).}
		If $M$ is closed and $p>1$, set
		\[
		I(t)=2\int_M\sum_j\frac{c_j}{m_j+1}u^{m_j+1}\,\dd V,
		\qquad
		D(t)=-\int_M\Bigl|\nabla\sum_jc_ju^{m_j}\Bigr|^p\,\dd V.
		\]
		Then $D/I^\gamma$ is nondecreasing on $(0,T)$ for every
		$0<\gamma\le pm_1/(m_1+1)$.
	\end{enumerate}
\end{corollary}

The simplest case consists of two powers. Let $m,\gamma>0$ with
$m\ne\gamma$, let $c>0$, and consider
\begin{equation}\label{ms:examples:two-powers}
	u_t=\Delta_p\bigl(u^m+cu^\gamma\bigr).
\end{equation}
Equations of this type arise in population dynamics: for $p=2$,
$m=2$, $\gamma=1$, \eqref{ms:examples:two-powers} is, after a rescaling
of time, the single-species dispersal law
$u_t=\Delta[(\alpha+\beta u)u]$ of Shigesada--Kawasaki--Teramoto
\cite{SKT}, with $c=\alpha/\beta$, and its limit $c\to0$ is the
porous-medium model $u_t=\Delta(u^2)$ of
\cite{GurneyNisbet,GurtinMacCamy}. In general, as $c\to0$,
\eqref{ms:examples:two-powers} becomes $u_t=\Delta_p(u^m)$; for $p=2$,
this is the porous-medium equation when $m>1$ and the fast-diffusion
equation when $0<m<1$. Put $\delta_1=m-1/(p-1)$ and $\delta_2=\gamma-1/(p-1)$. By
\eqref{ms:examples:sum-pressure},
\begin{equation}\label{ms:examples:two-b}
	\psi'=ms^{\delta_1}+c\gamma s^{\delta_2},
	\qquad
	b(s)=\delta_1+\frac{c\gamma(\delta_2-\delta_1)s^{\delta_2}}
	{ms^{\delta_1}+c\gamma s^{\delta_2}},
\end{equation}
and, when $\delta_1,\delta_2>0$, \eqref{ms:examples:sum-energies} and
\eqref{ms:examples:mixed-defect} read
\begin{equation}\label{ms:examples:two-energies}
	\begin{gathered}
		U_1(s)=\frac{m}{\delta_1(\delta_1+1)}s^{\delta_1+1},
		\qquad
		U_2(s)=\frac{c\gamma}{\delta_2(\delta_2+1)}s^{\delta_2+1},\\
		(U_1+U_2)D_{a_*}
		=(\delta_1-\delta_2)^2U_1U_2
		+\frac{\delta_1U_1+\delta_2U_2}{n}\cdot
		\frac{(\delta_*-\delta_1)U_1+(\delta_*-\delta_2)U_2}{\delta_*}>0,
	\end{gathered}
\end{equation}
where $\delta_*=\max\{\delta_1,\delta_2\}$ and
$a_*=n\delta_*/(n(p-1)\delta_*+p)$.

\begin{corollary}[Two powers]\label{ms:cor:two-powers}
	Let $u>0$ be a smooth solution of \eqref{ms:examples:two-powers} with
	values in $\mathcal U$, and let $b$ be given by
	\eqref{ms:examples:two-b}.
	\begin{enumerate}[label=\rm(\roman*)]
		\item \emph{Li--Yau estimate.}
		If $p\ge2$, $M$ is closed with $\Ric\ge0$, and $b(u_-)>-1/n$, then
		\[
		\frac{|\nabla(u^m+cu^\gamma)|^p}{u^{p/(p-1)}}
		-p\,\frac{mu^{m-1}+c\gamma u^{\gamma-1}}{u^{1/(p-1)}}\,u_t
		\le\frac{p\max_{s\in\mathcal U}\bigl(ms^{\delta_1}+c\gamma s^{\delta_2}\bigr)}
		{\bigl(b(u_-)+1/n\bigr)\,t}
		\qquad\text{on }M\times(0,T).
		\]
		\item \emph{Hamilton-type estimate.}
		Under the hypotheses of Corollary~\ref{ms:cor:sums}{\rm(ii)}, one has
		$u^{-1/(p-1)}|\nabla(u^m+cu^\gamma)|(\alpha-v)^{-q/p}
		\le C(R^{-1}+T^{-1/p}+\sqrt\kappa)$ on $Q_{R/2,T/2}$.
		\item \emph{$\mathcal W$-entropy formula.}
		Suppose that $m,\gamma>1/(p-1)$, that $a\ge a_*$, that $M$ is closed
		with $\Ric\ge0$, and that $u$ satisfies the regularity hypothesis of
		Theorem~\ref{ms:thm:entropy}. With $U_1,U_2$ as in
		\eqref{ms:examples:two-energies}, set
		\[
		D_a=\delta_1^2U_1+\delta_2^2U_2+\frac{\delta_1U_1+\delta_2U_2}{n}
		-\frac{a+1}{pa}\frac{(\delta_1U_1+\delta_2U_2)^2}{U_1+U_2}.
		\]
		Then
		\[
		\mathcal W_{F,p,a}(t)=t^{a+1}\int_M\left[
		\frac{|\nabla(u^m+cu^\gamma)|^p}{u^{1/(p-1)}}
		-\frac{a+1}{t}\bigl(U_1(u)+U_2(u)\bigr)\right]\dd V
		\]
		is strictly decreasing on $(0,T)$, and
		\[
		\begin{aligned}
			\frac{\dd}{\dd t}\mathcal W_{F,p,a}
			={}&-pt^{a+1}\int_M(\delta_1U_1+\delta_2U_2)
			\Bigl(\left||\nabla v|^{p-2}\nabla_i\nabla_jv
			+\frac{a(U_1+U_2)}{nt(\delta_1U_1+\delta_2U_2)}a_{ij}\right|_A^2\\
			&\qquad+|\nabla v|^{2p-4}\Ric(\nabla v,\nabla v)\Bigr)\dd V\\
			&-pt^{a+1}\int_M\bigl(\delta_1^2U_1+\delta_2^2U_2-D_a\bigr)
			\left(\Delta_pv+\frac{a(U_1+U_2)}{t(\delta_1U_1+\delta_2U_2)}\right)^2\dd V\\
			&-pt^{a+1}\int_MD_a(\Delta_pv)^2\,\dd V,
		\end{aligned}
		\]
		where $U_1,U_2,D_a$ are evaluated at $u$.
		\item \emph{Entropy power.}
		Under the hypotheses of {\rm(iii)},
		$\bigl(\int_M(U_1(u)+U_2(u))\,\dd V\bigr)^{-1/a}$ is nondecreasing and
		concave on $(0,T)$.
		\item \emph{Frequency monotonicity.}
		If $M$ is closed and $p>1$, then
		\[
		-\int_M|\nabla(u^m+cu^\gamma)|^p\,\dd V\cdot
		\left[2\int_M\left(\frac{u^{m+1}}{m+1}
		+\frac{cu^{\gamma+1}}{\gamma+1}\right)\dd V\right]^{-\lambda}
		\]
		is nondecreasing on $(0,T)$ for every
		$0<\lambda\le p\min\{m,\gamma\}/(\min\{m,\gamma\}+1)$.
	\end{enumerate}
\end{corollary}

\begin{remark}\label{ms:rem:two-powers}
	As $c\to0$, $b\to\delta_1$ and $\psi'\to ms^{\delta_1}$ uniformly on
	$\mathcal U$, and part~(i) becomes the Li--Yau estimate for
	$u_t=\Delta_p(u^m)$ under $\delta_1>-1/n$; for $p=2$, this covers the
	porous-medium equation and the fast-diffusion equation with
	$1-1/n<m<1$. If $1/(p-1)<\gamma<m$, then $a_*$ is the Wang--Chen
	exponent in \eqref{ms:intro:WC-entropy}, and
	$D_{a_*},U_2=O(c)$ by \eqref{ms:examples:two-energies}; thus
	part~(iii) with $a=a_*$ tends to \eqref{ms:intro:WC-formula}
	(Remark~\ref{ms:rem:entropy-models}), and for $c>0$ the term
	$cu^\gamma$ adds the strictly positive defect $D_{a_*}$. If
	$\gamma>m$, the exponent in (v) is the Sun--S\"urig exponent
	$pm/(m+1)$.
	For the dispersal law of \cite{SKT} ($p=2$, $m=2$,
	$\gamma=1$), one has $b(s)=2s/(2s+c)\in[0,1)$. Hence part~(i) holds on
	every density interval, with right-hand side
	$2(2u_++c)/\bigl((b(u_-)+1/n)t\bigr)$, and part~(v) holds with
	$\lambda=1$, so that $\log\int_M(u^3/3+cu^2/2)\,\dd V$ is convex in
	$t$. Parts~(iii) and (iv) require $m,\gamma>1/(p-1)$ and therefore
	cover neither this law nor the fast-diffusion range.
\end{remark}

\subsection{Arctangential diffusion}
\label{ms:examples:arctan}

Let $F(s)=\arctan s$. For $p=2$,
\[
\partial_tu=\Delta\arctan u
\]
is the arctangential heat equation of Brenier \cite{Brenier}, related
to the minimal surface equation in Minkowski space. Its diffusivity
$(1+u^2)^{-1}$ is close to $1$ at small densities, where the equation
behaves like the heat equation, and decays like $u^{-2}$ at large
densities.

Here
\[
v'(s)=\frac{s^{-1/(p-1)}}{1+s^2},
\qquad
\psi'=\frac{s^{(p-2)/(p-1)}}{1+s^2},
\qquad
b=\frac{p-2}{p-1}-\frac{2s^2}{1+s^2},
\qquad
b'(s)=-\frac{4s}{(1+s^2)^2}<0,
\]
so $b$ decreases from $(p-2)/(p-1)$ at $s=0$ to $(p-2)/(p-1)-2$ as
$s\to\infty$, and $\min_{\mathcal U}b=b(u_+)$. We now check the
hypotheses of the main theorems.

\emph{Gradient estimates.} By \eqref{ms:examples:LY-explicit}, the
Li--Yau estimate requires $b(u_+)>-1/n$, and Theorem~\ref{ms:thm:AB}
applies with $\varrho=(p-1)b(u_+)+p/n$. These conditions are
equivalent to
\begin{equation}\label{ms:examples:arctan-gradient-ranges}
	u_+^2<\frac{n(p-2)+p-1}{p(n-1)+1}
	\quad\text{(Li--Yau)},
	\qquad
	n=1\ \ \text{or}\ \ u_+^2<\frac{n(p-2)+p}{p(n-1)}
	\quad\text{(Aronson--B\'enilan)}.
\end{equation}
The Hamilton condition holds on short density intervals, as shown in
Section~\ref{ms:examples:mixed-powers}.

\emph{Entropy.} Let $p>2$. Then the pressure integral converges at
zero, and with $v(0)=0$,
\begin{equation}\label{ms:examples:arctan-energies}
	\begin{gathered}
		v(s)=\int_0^s\frac{r^{-1/(p-1)}}{1+r^2}\,\dd r,
		\qquad
		U(s)=\int_0^s\frac{(s-r)r^{-1/(p-1)}}{1+r^2}\,\dd r,
		\qquad
		P(s)=\int_0^s\frac{r^{(p-2)/(p-1)}}{1+r^2}\,\dd r,\\
		D_a(s)=\frac{s^{(2p-3)/(p-1)}}{1+s^2}
		-\left(1-\frac1n\right)P(s)
		-\frac{a+1}{pa}\frac{P(s)^2}{U(s)}.
	\end{gathered}
\end{equation}
As $s\downarrow0$, $U(s)\sim\frac{(p-1)^2}{(p-2)(2p-3)}s^{(2p-3)/(p-1)}$
and $P(s)\sim\frac{p-1}{2p-3}s^{(2p-3)/(p-1)}$, so that
\[
\lim_{s\downarrow0}\frac{D_a(s)}{P(s)}
=\frac{p-2}{p-1}+\frac1n-\frac{(a+1)(p-2)}{pa(p-1)},
\]
which is positive if and only if
$a>a_0:=\dfrac{n(p-2)}{(p-1)[n(p-2)+p]}$. For such $a$,
\[
s_0=s_0(n,p,a):=\sup\{s>0:D_a>0\text{ on }(0,s]\}\in(0,\infty].
\]
If $n\ge2$, then $s_0<\infty$: as $s\to\infty$, $P(s)$ has a finite
limit and $sP'(s)=s^{(2p-3)/(p-1)}/(1+s^2)\to0$, so that
$D_a\le sP'-(1-\frac1n)P<0$ for large $s$.

\emph{Frequency.} Let $\ell\in\mathcal U$. Since
$(1+u_+^2)^{-1}\le F'\le(1+u_-^2)^{-1}$ on $\mathcal U$, integrating
$\frac{\dd}{\dd s}[F(s)-F(\ell)]^2=2[F(s)-F(\ell)]F'(s)$ from $\ell$
gives $[F(s)-F(\ell)]^2\le2(1+u_-^2)^{-1}\Psi_\ell(s)$, so
\eqref{ms:freq:structure} holds with $\kappa=2(1+u_+^2)/(1+u_-^2)$, where
\[
\Psi_\ell(s)=s(\arctan s-\arctan\ell)-\frac12\log\frac{1+s^2}{1+\ell^2}.
\]

\begin{corollary}[Arctangential diffusion]\label{ms:cor:arctan}
	Let $F(s)=\arctan s$, and let $u>0$ be a smooth solution of
	\eqref{ms:intro:PF} with values in $\mathcal U$.
	\begin{enumerate}[label=\rm(\roman*)]
		\item \emph{Li--Yau estimate (Theorem~\ref{ms:thm:LY}).}
		If $p\ge2$, $M$ is closed with $\Ric\ge0$, and the first condition in
		\eqref{ms:examples:arctan-gradient-ranges} holds, then, on
		$M\times(0,T)$,
		\[
		\frac{|\nabla u|^p}{u^{p/(p-1)}(1+u^2)^p}
		-\frac{p\,u_t}{u^{1/(p-1)}(1+u^2)}
		\le\frac{p\max_{s\in\mathcal U}s^{(p-2)/(p-1)}(1+s^2)^{-1}}
		{\Bigl(\frac{p-2}{p-1}-\frac{2u_+^2}{1+u_+^2}+\frac1n\Bigr)t}.
		\]
		\item \emph{Aronson--B\'enilan estimate (Theorem~\ref{ms:thm:AB}).}
		If $p\ge2$, $M$ is closed with $\Ric\ge0$, and the second condition
		in \eqref{ms:examples:arctan-gradient-ranges} holds, then
		$\varrho:=p-2+\frac pn-\frac{2(p-1)u_+^2}{1+u_+^2}>0$ and
		\[
		\Delta_pv\ge-\frac1{\varrho t},
		\qquad\text{equivalently}\qquad
		\frac{u_t}{u}\ge\frac{|\nabla u|^p}{u^2(1+u^2)^{p-1}}-\frac1{\varrho t}
		\qquad\text{on }M\times(0,T).
		\]
		\item \emph{Hamilton-type estimate (Theorem~\ref{ms:thm:Hamilton}).}
		Under the hypotheses of Corollary~\ref{ms:cor:sums}{\rm(ii)}, with
		$\eta,\varepsilon_0$ determined by $F=\arctan$ on $[u_0/2,2u_0]$,
		\[
		\sup_{Q_{R/2,T/2}}
		\frac{|\nabla u|}{u^{1/(p-1)}(1+u^2)(\alpha-v)^{q/p}}
		\le C\left(\frac1R+\frac1{T^{1/p}}+\sqrt\kappa\right),
		\]
		where $v$ is the pressure \eqref{ms:intro:pressure} and $C$ depends
		only on $n,p,q,\varepsilon$ and $u_0$.
		\item \emph{$\mathcal W$-entropy formula
			(Theorem~\ref{ms:thm:entropy}).}
		Let $p>2$, $a>a_0$, and $u_+\le s_0(n,p,a)$. Suppose that $M$ is
		closed with $\Ric\ge0$ and that $u$ satisfies the regularity
		hypothesis of Theorem~\ref{ms:thm:entropy}. With $U,P,D_a$ as in
		\eqref{ms:examples:arctan-energies},
		\[
		\mathcal W_{F,p,a}(t)
		=t^{a+1}\int_M\left[
		\frac{|\nabla u|^p}{u^{1/(p-1)}(1+u^2)^p}
		-\frac{a+1}{t}U(u)\right]\dd V
		\]
		is strictly decreasing on $(0,T)$, and
		\[
		\begin{aligned}
			\frac{\dd}{\dd t}\mathcal W_{F,p,a}
			={}&-pt^{a+1}\int_MP(u)
			\left(\left||\nabla v|^{p-2}\nabla_i\nabla_jv
			+\frac{aU(u)}{ntP(u)}a_{ij}\right|_A^2
			+|\nabla v|^{2p-4}\Ric(\nabla v,\nabla v)\right)\dd V\\
			&-pt^{a+1}\int_M\left[\frac{u^{(2p-3)/(p-1)}}{1+u^2}-P(u)-D_a(u)\right]
			\left(\Delta_pv+\frac{aU(u)}{tP(u)}\right)^2\dd V\\
			&-pt^{a+1}\int_MD_a(u)(\Delta_pv)^2\,\dd V.
		\end{aligned}
		\]
		\item \emph{Entropy-power concavity
			(Theorem~\ref{ms:thm:entropy-power}).}
		Under the hypotheses of {\rm(iv)},
		$\mathcal N_{F,p,a}(t)=\bigl(\int_MU(u)\,\dd V\bigr)^{-1/a}$ is
		nondecreasing and concave on $(0,T)$, and
		$\mathcal N_{F,p,a}''(t_0)=0$ if and only if $u(\cdot,t_0)$ is
		constant on each connected component of $M$.
		\item \emph{Frequency monotonicity (Theorem~\ref{ms:thm:frequency}).}
		If $M$ is closed, $p>1$ and $\ell\in\mathcal U$, set
		\[
		I(t)=2\int_M\left[u(\arctan u-\arctan\ell)
		-\frac12\log\frac{1+u^2}{1+\ell^2}\right]\dd V,
		\qquad
		D(t)=-\int_M\frac{|\nabla u|^p}{(1+u^2)^p}\,\dd V.
		\]
		Then $D/I^\gamma$ is nondecreasing on every time interval on which
		$I>0$, for every $0<\gamma\le p(1+u_-^2)/[2(1+u_+^2)]$.
	\end{enumerate}
\end{corollary}

\begin{remark}
	For $p=2$ and $n\ge2$, \eqref{ms:examples:arctan-gradient-ranges}
	reads $u_+^2<1/(2n-1)$ for the Li--Yau estimate and $u_+^2<1/(n-1)$
	for the Aronson--B\'enilan estimate; thus the monotonicity of $b$ gives
	a larger range. The restriction $u_+\le s_0$ in part~(iv) is
	structural: for $n\ge2$, McCann's quantity $sP'-(1-\frac1n)P$ is
	negative at large densities, so the internal energy of the
	arctangential law is displacement convex in dimension $n$ only on
	bounded density ranges.
\end{remark}

\section*{AI tool disclosure}
	During the preparation of this manuscript, the authors
	used AI tools, primarily ChatGPT-6 Astra, for \LaTeX\
	editing, improving the presentation of proofs, and
	checking calculations and consistency.
	The research ideas, mathematical approach, and main
	proof arguments were developed by the authors.
	The authors take full responsibility for the results,
	references, and final manuscript.

\begingroup
\small\linespread{1}\selectfont
Yu-Zhao Wang\\
\emph{E-mail:} wangyuzhao@sxu.edu.cn\\
School of Mathematics and  Statistics, Shanxi University, Taiyuan, 030006, Shanxi, China 

Jian-Hua Hao \\
\emph{E-mail:} haojianhua@sxu.edu.cn\\
School of Mathematics and  Statistics, Shanxi University, Taiyuan, 030006, Shanxi, China.
\par\endgroup

\end{document}